\documentclass[final,3p]{elsarticle}

\usepackage{siunitx}
\usepackage{graphicx}
\usepackage{amssymb}
\usepackage{caption}
\usepackage{mathtools}
\usepackage{tikz}
\usetikzlibrary{positioning}
\usepackage{amsthm}
\usepackage{amsfonts}
\usepackage{amsmath}
\usepackage{BOONDOX-cal}
\usepackage{xfrac}
\usepackage{relsize}
\usepackage{float}
\usepackage{fixltx2e}
\usepackage[justification=centering]{caption}
\usepackage{hyphenat}
\usepackage{epsfig}
\usepackage{epstopdf}
\usepackage{lineno}
\usepackage{subcaption}
\usepackage{hhline}
\usepackage{url}
\modulolinenumbers[5]
\usepackage{changes}
\usepackage{soul}
\usepackage{xcolor}
\usepackage{hyperref}
\usepackage{setspace} 
\hypersetup{
    colorlinks=true,
    linkcolor=blue,
    filecolor=magenta,      
    urlcolor=cyan,
}

\usepackage{algorithm}
\usepackage{algpseudocode}

\DeclareMathOperator*{\argmin}{arg\,min}

\journal{Arxiv}
\begin{document}

\begin{frontmatter}

\author[rvt]{Stefanos Nikolopoulos\corref{cor1}}
\ead{stefnikolopoulos@mail.ntua.gr}

\author[rvt]{Ioannis Kalogeris}
\ead{ikalog@mail.ntua.gr}

\author[rvt]{George Stavroulakis}
\ead{stavroulakis@nessos.gr}

\author[rvt]{Vissarion Papadopoulos}
\ead{vpapado@central.ntua.gr}
\ead[url]{mgroup.ntua.gr}

\cortext[cor1]{corresponding author}
\address[rvt]{MGroup, Engineering Simulations Lab, Institute of Structural Analysis and Antiseismic research, National Technical University of Athens, Zografou Campus, 9 Iroon Polytechniou str, 15780 Zografou}

\title{AI-enhanced iterative solvers for accelerating the solution of large-scale parametrized systems}

\begin{abstract}
Recent advances in the field of machine learning open a new era in high performance computing for challenging computational science and engineering applications. In this framework, the use of advanced machine learning algorithms for the development of accurate and cost-efficient surrogate models of complex physical processes has already attracted major attention from scientists. However, despite their powerful approximation capabilities, surrogate model predictions are still far from being near to the `exact' solution of the problem. To address this issue, the present work proposes the use of up-to-date machine learning tools in order to equip a new generation of iterative solvers of linear equation systems, capable of very efficiently solving large-scale parametrized problems at any desired level of accuracy. The proposed approach consists of the following two steps. At first, a reduced set of model evaluations is performed using a standard finite element methodology and the corresponding solutions are used to establish an approximate mapping from the problem's parametric space to its solution space using a combination of deep feedforward neural networks and convolutional autoencoders. This mapping serves a means to obtain very accurate initial predictions of the system's response to new query points at negligible computational cost. Subsequently, an iterative solver inspired by the Algebraic Multigrid method in combination with Proper Orthogonal Decomposition, termed POD-2G, is developed that successively refines the initial predictions of the surrogate model towards the exact solution. The application of POD-2G as a standalone solver or as preconditioner in the context of preconditioned conjugate gradient methods is demonstrated on several numerical examples of large scale systems, with the results indicating its strong superiority over conventional iterative solution schemes.


\end{abstract}

    \end{frontmatter}
    
\section{Introduction}

In scientific computing, there is a constant need for solving larger and computationally more demanding problems with increased accuracy and improved numerical performance. This holds particularly true in multi-query scenarios such as optimization, uncertainty quantification, inverse problems and optimal control, where the problems under investigation need to be solved for numerous different parameter instances with high accuracy and efficiency. In this regard, constructing efficient numerical solvers for complex systems described by partial differential equations is crucial for many scientific disciplines. The preconditioned conjugate gradient method (PCG) \cite{Barrett1994,Benzi2000, Lin2014,Herzog2010} and the preconditioned generalised minimal residual method (PGMRES) \cite{Saad1986,Shakib1989, Baglama1998} are amongst the most powerful and versatile approaches to treat such problems. In these methods, the choice of a suitable preconditioner plays a major role on the convergence and scalability of the solvers and notable examples include the incomplete Choleski factorization \cite{Concus1985} and domain decomposition methods \cite{Toselli2005, TALLEC1991}, such as the popular FETI methods \cite{Farhat1991, FARHAT1994, Fragakis2003} and the additive Schwarz methods \cite{Cai1999, DAAS2021}. In a similar fashion, Algebraic and Geometric Multigrid  (AMG, GMG, resp.) \cite{Trottenberg2000} are equally well-established methods that are commonly employed for accelerating standard iterative solvers and may also service as highly efficient preconditioners for PCG \cite{Iwamura2003, Heys2005,Langer2003} or PGMRES \cite{RAMAGE1999, Wienands2000, Vakili2009}. 

Nevertheless, optimizing the aforementioned solvers so as to attain a uniformly fast convergence for multiple parameter instances, as required in multi-query problems, remains a challenging task to this day. To tackle this problem, several works suggest the use of interpolation methods tasked with constructing approximations of the system's inverse operator for different parameter values \cite{Zahm2016,Bergamaschi2020, Carr2021}, which can then be used as preconditioners. Another approach can be found in \cite{STAVROULAKIS2014}, where primal and dual FETI decomposition methods with customized preconditioners are developed in order to accelerate the solution of stochastic problems in the context of Monte Carlo simulation, as well as intrusive Galerkin methods. Augmented Krylov Subspace methods showed great promise in handling sequences of linear systems \cite{Saad1997}, such as those arising in parametrized PDEs, however, the augmentation of the usual Krylov subspace with data from multiple previous solves led in certain cases to disproportional computational and memory requirements. To alleviate this cost, optimal truncation strategies have been proposed in \cite{Sturler1999}, as well as deflation techniques \cite{Chapman1997,Saad2000, Daas2021a}.   

In recent days, the rapid advancements in the field of machine learning (ML) have offered researchers new tools to tackle challenging problems in multi-query scenarios. For instance, deep feedforward neural networks (FFNNs) have been successfully employed to construct response surfaces of quantities of interest in complex problems \cite{Papadrakakis1996,PAPADRAKAKIS2002, SEYHAN2005, HOSNIELHEWY2006,Chojaczyk2015}. Convolutional neural networks (CNNs) in conjuction with FFNNs have been employed to predict the high-dimensional system response at different parameter instances \cite{NIKOLOPOULOS2021, NIKOLOPOULOS2022,XU2020}. In addition, recurrent neural networks demonstrated great potential in transient problems for propagating the state of the system forward in time without the need of solving systems of equations \cite{Yu2019,Zhou2019}. All these non-intrusive approaches utilize a reduced set of system responses to build an emulator of the system's input-output relation for different parameter values. As such, they are particularly cheap to evaluate and can be very accurate in certain cases. However, these methods can be characterized as physics-agnostic in the sense that the derived solutions do not satisfy any physical laws. This problem is remedied to some extent from intrusive approaches based on reduced basis methods, such as Principal Orthogonal Decomposition (POD) \cite{Carlberg2011,Zahr2017,Agathos2020} and proper Generalized Decomposition \cite{Chinesta2010,Ladeveze2010,Ladeveze2011}. These methods rely on the premise that a small set of appropriately selected basis vectors suffices to construct a low-dimensional subspace of the system's high-dimensional solution space and the projection of the governing equations to this subspace will come at minimum error. In addition, several recent works have investigated the combination of either linear or nonlinear dimensionality reduction algorithms and non-intrusive interpolation schemes to construct cheap emulators of complex systems \cite{DALSANTO2020, SALVADOR2021,KALOGERIS2021, dosSantos2022, Kadeethum2022,VLACHAS2021, Heaney, LAZZARA2022107629}. Nevertheless, none of these surrogate modelling schemes can guarantee convergence to the exact solution of the problem.

In the effort to combine the best of two worlds, a newly emergent research direction is that of enhancing linear algebra solvers with machine-learning algorithms. For instance, POD has been successfully employed to truncate the augmented Krylov subspace and retain only the high-energy modes \cite{Carlberg2016} for efficiently solving sequences of linear systems of equations characterized by varying right-hand sides and symmetric-positive-definite matrices. In \cite{Heinlein2019}, neural networks were trained for predicting the geometric location of constraints in the context of domain decomposition methods, leading to enhanced algorithm robustness. Moreover, the close connection between multigrid methods and CNNs has been studied in several recent works, which managed to accelerate their convergence by providing data-driven smoothers \cite{CHEN2022}, prolongation and restriction operators \cite{luz2020}.

The present work aims at bridging the gap between machine learning and linear algebra algorithms for accelerating the solution of real-life computational mechanics problems in multi-query scenarios. To this end, a novel strategy is proposed to utilize ML tools in order to obtain system solutions within a prescribed accuracy threshold, with faster convergence rates than conventional solvers. The proposed approach consists of two steps. Initially, a reduced set of model evaluations is performed and the corresponding solutions are used to establish an approximate mapping from the problem's parametric space to its solution space using a combination of deep FFNNs and CAEs. This mapping serves a means of acquiring very accurate initial predictions of the system's response to new query points at negligible computational cost. The error in these predictions, however, may or may not satisfy the prescribed accuracy threshold. Therefore, a second step is proposed herein, which further utilizes the knowledge from the already available system solutions, in order to construct a data-driven iterative solver. This solver is inspired by the idea of the Algebraic Multigrid method combined with Proper Orthogonal Decomposition, termed POD-2G, that successively refines the initial prediction of the surrogate model towards the exact system solutions with significantly faster convergence rates.

The paper is organised as follows. In Section \ref{sec2} the basic principles of the PCG and AMG iterative solvers are illustrated. In Section \ref{sec3}, the elaborated methodology for developing an AI-enhanced linear algebra solver is presented. Section \ref{sec4} presents a series of numerical examples that showcase the performance of the method compared to conventional iterative solvers. Section \ref{sec5} summarizes the outcomes of this work and discusses possible extensions.

\section{Iterative solvers for FE systems} \label{sec2}

\subsection{The Finite Element Method} \label{sec2.0}
This work focuses on linear elliptic PDEs defined on a domain $\Omega \subseteq \mathbb{R}^{dim}$, $dim=1,2,3$, which are parametrized by a vector of parameters $\boldsymbol{\theta}\in \boldsymbol{\Theta}$, with $\boldsymbol{\Theta}\subseteq\mathbb{R}^n$ being the parameter space. The variational formulation of the PDE can be stated as: given $\boldsymbol{\theta}\in \boldsymbol{\Theta}$, find the solution $v=v(\boldsymbol{\theta})$ from the Hilbert space $\mathcal{V}=\mathcal{V}(\Omega)$ such that

\begin{equation} \label{eq:weakForm}
    \kappa\left(v ,w; \boldsymbol{\theta}\right)=f \left(w;\boldsymbol{\theta}\right)
\end{equation}

\noindent for every $w\in \mathcal{V}(\Omega)$ with compact support in $\Omega$. The Lax-Milgram lemma proves that eq. \eqref{eq:weakForm} has a unique solution for every  $\boldsymbol{\theta}$, provided that the bilinear form $\kappa(\cdot,\cdot;\boldsymbol{\theta})$ is continuous and coercive and $f\left(\cdot;\boldsymbol{\theta} \right)$ is a continuous one-form. 

In practice, however, obtaining an exact solution $v$ is not feasible for most applications of interest and instead, an approximate solution is sought using numerical techniques, such the finite element method (FEM). In FEM, a finite-dimensional subspace $\mathcal{V}_h \subseteq \mathcal{V}$ is considered, which is spanned by a finite number of polynomial basis vectors $\lbrace N_i  \rbrace_{i=1}^{\bar{N}}$. These polynomials are compactly supported on a set of small polyhedra (finite elements) that partition the domain $\Omega$ and within each element $e$ the approximate displacement vector field $v_h\in\mathcal{V}_h$ and test functions $w_h$ are expressed as:

\begin{align}
    v_h^e &=\sum_{i=1}^{\bar{N}} \mathrm{u}_i^e N_i^e  \\
    w_h^e &=\sum_{i=1}^{\bar{N}} \mathrm{w}_i^e N_i^e 
\end{align}

\noindent where $\boldsymbol{u}^e=[\mathrm{u}_i^e,\cdots,\mathrm{u}_{\bar{N}}^e]^T\in \mathbb{R}^{\bar{N}}$ are the coefficients in the expansion of the unknown field approximation, obtained using a Galerkin minimization that relies on the linearity of the forms $\kappa,f$ and the orthogonality of the polynomial basis vectors. Since eq. \eqref{eq:weakForm} must hold within each finite element $e$ and for any test function $w$, the system of linear equations follows:

\begin{equation}
    \kappa\left(\sum_{j=1}^{\bar{N}} \mathrm{u}_j^e N_j, N_i; \boldsymbol{\theta}\right)=f \left(N_i;\boldsymbol{\theta}\right), \ \ \text{for} \ i=1,...,\bar{N}
\end{equation}

or, due to the linearity of $\kappa$,

\begin{equation}\label{linearSystemInElement}
    \sum_{j=1}^{\bar{N}}\kappa\left( N_j, N_i; \boldsymbol{\theta}\right)  \mathrm{u}_j^e=f \left(N_i;\boldsymbol{\theta}\right), \ \ \text{for} \ i=1,...,\bar{N}
\end{equation}

Equation \eqref{linearSystemInElement} describes an $\bar{N} \times \bar{N}$ linear system of equations to be satisfied within the $e$-th element. Repeating this procedure for all elements and appropriately assembling the respective equations will result in the following $d\times d$ linear system

\begin{equation} \label{eq:fineProblem}
    \boldsymbol{K}(\boldsymbol{\theta})\boldsymbol{u}(\boldsymbol{\theta})=\boldsymbol{f}(\boldsymbol{\theta})
\end{equation}

\noindent with $d$ being the total number of unknowns in the system, $\boldsymbol{K}\in\mathbb{R}^{d\times d}$ is a real symmetric positive definite matrix, $\boldsymbol{u}\in\mathbb{R}^d$ is the unknown solution vector and $\boldsymbol{f}\in\mathbb{R}^d$ the force vector. 

Solving such a linear system for a detailed discretization ($d \gg 1$) can be computationally intensive, particularly in multiquery problems that require numerous system evaluations for various instances of parameters $\boldsymbol{\theta}$, such as optimization, parameter inference, uncertainty propagation, sensitivity analysis, etc. Therefore, it becomes evident that efficient numerical solvers for linear systems of equations are of vital importance in the analysis of large scale real-world problems. The following section revisits the basic ideas behind two of most efficient methods for solving such systems, namely, the PCG and the AMG methods.

\subsection{Preconditioned conjugate gradient method} \label{sec2.1}
The Conjugate Gradient method was originally proposed by Hestenes and Stiefel as a direct method \cite{Hestenes1952} for solving linear systems, but its full potential was demonstrated in the frame of iterative solvers for large-scale sparse systems of the form $\boldsymbol{K}\boldsymbol{u}=\boldsymbol{f}$, with $\boldsymbol{K}$ being a symmetric positive definite matrix. The goal of CG is to minimize the quadratic function

\begin{equation}
    Q(\boldsymbol{u})=\frac{1}{2}\boldsymbol{u}^T\boldsymbol{K}\boldsymbol{u}-\boldsymbol{f}^T\boldsymbol{u}
\end{equation}

\noindent which is equivalent to setting the residual $\boldsymbol{r}=-\nabla Q(\boldsymbol{u})=\boldsymbol{f}-\boldsymbol{K}\boldsymbol{u}$ to zero. 

Let us assume an initial guess $\boldsymbol{u}^{(0)}$ for the system, which, in the absence of any other information, is taken $\boldsymbol{u}^{(0)}=\boldsymbol{0}$, with corresponding residual $\boldsymbol{r}^{(0)}=\boldsymbol{f}$.  Then, we can consider the Krylov subspaces,

\begin{equation}
    \mathcal{K}_0=\lbrace \boldsymbol{0} \rbrace, \ \ \ \mathcal{K}_k=span\lbrace \boldsymbol{f},\boldsymbol{K}\boldsymbol{f},\dots,\boldsymbol{K}^{k-1}\boldsymbol{f} \rbrace, \ \ \text{for } k\geq 1
\end{equation}

\noindent These subspaces are nested, $\mathcal{K}_0\subseteq\mathcal{K}_1\subseteq \dots$, and have the key property that $\boldsymbol{K}^{-1}\boldsymbol{f} \in \mathcal{K}_d$. Then, a Krylov sequence $\lbrace \boldsymbol{u}^{(1)},\boldsymbol{u}^{(2)},\cdots \rbrace$ consists of the vectors $\boldsymbol{u}^{(k)}$ such that

\begin{equation} \label{eq:MinProblemCG}
    \boldsymbol{u}^{(k)}=\argmin_{\boldsymbol{u}\in\mathcal{K}_k} Q(\boldsymbol{u}), \ \ k=1,2,\dots
\end{equation}

\noindent Bbased on the previous property, it follows that $\boldsymbol{u}^{(d)}=\boldsymbol{K}^{-1}\boldsymbol{f}$. In this regard, CG is a recursive method for computing the Krylov sequence $\lbrace \boldsymbol{u}^{(0)},\boldsymbol{u}^{(1)},\dots \rbrace$. It can be proven that the corresponding (nonzero) residuals $\boldsymbol{r}^{(k)}=\boldsymbol{f}-\boldsymbol{K}\boldsymbol{u}^{k()}$ form an orthogonal basis for the Krylov subspaces, that is

\begin{equation}
    \mathcal{K}_k=span\lbrace \boldsymbol{r}^{(0)},\boldsymbol{r}^{(1)},\dots,\boldsymbol{r}^{(k-1)} \rbrace, \ \ \ \left(\boldsymbol{r}^{(j)}\right)^T\boldsymbol{r}^{(i)}=0, \text{ for } i\neq j
\end{equation}

\noindent and a sequence of conjugate ($\boldsymbol{K}$-orthogonal) basis vectors $\boldsymbol{p}_k$ can be obtained by applying the Gram-Schmidt process to the $\boldsymbol{r}^{(k)}$ vectors as follows:

\begin{equation}
\boldsymbol{p}_0=\boldsymbol{r}^{(0)}, \ \ \ \boldsymbol{p}_{k}=\boldsymbol{r}^{(k)}-\sum_{i<k}\frac{\boldsymbol{p}_{i}^T\boldsymbol{K}\boldsymbol{r}^{(k)}}{\boldsymbol{p}_{i}^T\boldsymbol{K}\boldsymbol{p}_{i}}\boldsymbol{p}_{i}, \ k=1,2,\dots
\end{equation}

\noindent or, equivalently, 

\begin{align}
 \boldsymbol{p}_{k}&=\boldsymbol{r}^{(k)}-\frac{\boldsymbol{p}_{k-1}^T\boldsymbol{K}\boldsymbol{r}^{(k)}}{\boldsymbol{p}_{k-1}^T\boldsymbol{K}\boldsymbol{p}_{k-1}}\boldsymbol{p}_{k-1}, \ k=1,2,\dots   \nonumber \\
 &=\boldsymbol{r}^{(k)}+\frac{\left(\boldsymbol{r}^{(k)}\right)^T\boldsymbol{r}^{(k)}}{\left(\boldsymbol{r}^{(k-1)}\right)^T\boldsymbol{r}^{(k-1)}}\boldsymbol{p}_{k-1}, \ k=1,2,\dots
\end{align}

The solution $\boldsymbol{u}^{(k+1)}=\argmin_{\boldsymbol{u}\in\mathcal{K}_{k+1}} Q(\boldsymbol{u})$ of eq. \eqref{eq:MinProblemCG} can be expressed as a linear combination of the basis vectors $\lbrace\boldsymbol{p_0},\dots,\boldsymbol{p_k} \rbrace$

\begin{equation}
    \boldsymbol{u}^{(k+1)}= \sum_{i=0}^k \alpha_i\boldsymbol{p}_i
\end{equation}

\noindent with the coefficients $\alpha_i$ obtained from the Galerkin projections:

\begin{equation}
\alpha_i=\frac{\boldsymbol{p}_i^T\boldsymbol{r}^{(i)}}{\boldsymbol{p}_i^T\boldsymbol{K}\boldsymbol{p}_i}
\end{equation}

Using the fact that $\boldsymbol{u}^{(k)}=\sum_{i=0}^{k-1} \alpha_i\boldsymbol{p}_i $, then, the Krylov sequence and the corresponding residuals are given by the relations:

\begin{align}
    \boldsymbol{u}^{(k+1)} &= \boldsymbol{u}^{(k)}+\alpha_k\boldsymbol{p}_k \label{eq:CG_u} \\
    \boldsymbol{r}^{(k+1)} &= \boldsymbol{r}^{(k)}-\alpha_k\boldsymbol{K}\boldsymbol{p}_k \label{eq:CG_r}
\end{align}

In the above, we could consider an initial guess $\boldsymbol{u}^{(0)}\neq \boldsymbol{0}$ and solve the system $\boldsymbol{K}\bar{\boldsymbol{u}}=\boldsymbol{f}-\boldsymbol{K}\boldsymbol{u}^{(0)}$, with $\boldsymbol{u}=\bar{\boldsymbol{u}}+\boldsymbol{u}^{(0)}$. This is the same as initializing the CG algorithm with $\lbrace \boldsymbol{u}^{(0)}$, $\boldsymbol{r}^{(0)}=\boldsymbol{f}-\boldsymbol{K}\boldsymbol{u}^{(0)}\rbrace$ and updating this guess according to equations \eqref{eq:CG_u}-\eqref{eq:CG_r} for $k=1,2,\dots$, until $\boldsymbol{r}^{(k)}$ is suffiently small. In theory, CG terminates in at most $d$ steps, however, due to rounding errors it may take more than $d$ steps or even fail in practice. Also, the improvement in the approximations $\boldsymbol{u}^{(k)}$ is determined by the condition number $c(\boldsymbol{K})$ of the system matrix $\boldsymbol{K}$; the larger $c(\boldsymbol{K})$ is, the slower the improvement.

A standard approach to enhance the convergence of the CG method is through preconditioning (PCG), namely the application of a linear transformation to the system with a matrix $\boldsymbol{T}$, called the preconditioner, in order to reduce the condition number of the problem. Thus, the original system $\boldsymbol{K}\boldsymbol{u}-\boldsymbol{f}=0$ is replaced with $\boldsymbol{T}^{-1}\left(\boldsymbol{K}\boldsymbol{u}-\boldsymbol{f}\right)=0$, such that $c(\boldsymbol{T}^{-1}\boldsymbol{K})$ is smaller than $c(\boldsymbol{K})$. The steps of the PCG algorithm are presented in algorithm \ref{alg:PCGalgorithm}.

\begin{algorithm}
\setstretch{1.0}
\caption{PCG algorithm}\label{alg:PCGalgorithm}
\begin{algorithmic}[1]
\State \textbf{Input:} $\boldsymbol{K}\in\mathbb{R}^{d\times d}$, rhs $\boldsymbol{f}\in\mathbb{R}^d$, preconditioner $\boldsymbol{T}\in\mathbb{R}^{d \times d}$, residual tolerance $\delta$ and an initial approximation $\boldsymbol{u}^{(0)}$
\State set $k=0$, initial residual $\boldsymbol{r}^{(0)}=\boldsymbol{f}-\boldsymbol{K}\boldsymbol{u}^{(0)}$
\State $\boldsymbol{s}_0=\boldsymbol{T}^{-1}\boldsymbol{r}^{(0)}$
\State $\boldsymbol{p}_0=\boldsymbol{s}_0$
\While{$\Vert \boldsymbol{r}^{(k)}\Vert < \delta$ }
    \State $\alpha_k=\frac{\left(\boldsymbol{r}^{(k)}\right)^T\boldsymbol{s}_k}{\boldsymbol{p}_k^T\boldsymbol{K}\boldsymbol{p}_k}$
    \State $\boldsymbol{u}^{(k+1)}=\boldsymbol{u}^{(k)}+\alpha_k \boldsymbol{p}_k$
    \State $\boldsymbol{r}^{(k+1)}=\boldsymbol{r}^{(k)}-\alpha_k \boldsymbol{K}\boldsymbol{p}_k$
    \State $\boldsymbol{s}_{k+1}=\boldsymbol{T}^{-1}\boldsymbol{r}^{(k+1)}$
        \State $\beta_k=\frac{\left(\boldsymbol{r}^{(k+1)}\right)^T\boldsymbol{s}_{k+1}}{\left(\boldsymbol{r}^{(k)}\right)^T\boldsymbol{s}_{k}}$
    \State $\boldsymbol{p}_{k+1}=\boldsymbol{s}_{k+1}+\beta_k\boldsymbol{p}_k$
    \State $k=k+1$
\EndWhile
\end{algorithmic}
\end{algorithm}

The choice of the preconditioner $\boldsymbol{T}$ in PCG plays a crucial role in the fast convergence of the algorithm. Some generic choices include the Jacobi (diagonal) preconditioner $\boldsymbol{T}=diag(\boldsymbol{K})$ and the incomplete Cholesky factorization $\boldsymbol{T}=\boldsymbol{\hat{L}}\boldsymbol{\hat{L}}^T$, with $\boldsymbol{\hat{L}}$ being a sparse lower triangular matrix such that $\boldsymbol{K}\approx\boldsymbol{\hat{L}}\boldsymbol{\hat{L}}^T$. Another popular choice is the incomplete LU factorization $\boldsymbol{T}=\boldsymbol{\tilde{L}}\boldsymbol{\tilde{U}}$, with $\boldsymbol{\tilde{L}}$ being a lower unitriangular matrix and $\boldsymbol{\tilde{U}}$ an upper triangular, such that $\boldsymbol{K}\approx\boldsymbol{\tilde{L}}\boldsymbol{\tilde{U}}$. Moreover, multigrid methods such as the AMG, elaborated on the next section, apart from standalone iterative schemes, are also very effective as preconditioners to the CG method.

\subsection{Algebraic Multigrid Method} \label{sec2.2}
AMG was originally introduced in the 1980's \cite{Ruge1987} as an efficient numerical approach for solving large ill-conditioned sparse linear systems and eigenproblems. Its main difference from the (geometric) multigrid method lies only in the method of coarsening. While multigrid methods require knowledge of the mesh, AMG methods extract all the needed information from the system matrix. AMG methods have been successfully applied to numerous problems including PDEs, sparse Markov chains and problems involving graph Laplacians (e.g. \cite{Stuben2001, Brezina2001, Treister2010,Napov2016, Facca2021}). The key idea in AMG algorithms is to employ a hierarchy of progressively coarser approximations to the linear system under consideration in order to accelerate the convergence of classical simple and cheap iterative processes, such as the damped Jacobi or Gauss-Seidel. These methods, commonly referred to as relaxation or smoothing, are very efficient in eliminating the high-frequency error modes, but inefficient in resolving the low-energy modes. AMG overcomes this problem through the coarse-level correction, as elaborated below.

Let us consider the linear system of eq. \eqref{eq:fineProblem}, which describes the fine problem and let $\boldsymbol{u}^{(0)}$ be an initial solution to it. The two-level AMG defines a prolongation operator $ \boldsymbol{P}$, which is a full-column rank matrix in $\mathbb{R}^{d\times d_c}$, $d_c < d$ and a relaxation scheme such as the Gauss-Seidel (GS). Then, the two-level AMG algorithm consists in the steps shown in algorithm \ref{alg:AMGalgorithm}:

\begin{algorithm}[H]
\setstretch{1.0}
\caption{Two-level AMG algorithm}\label{alg:AMGalgorithm}
\begin{algorithmic}[1]
\State \textbf{Input:} $\boldsymbol{K}\in\mathbb{R}^{d\times d}$, rhs $\boldsymbol{f}\in\mathbb{R}^d$, prolongation operator $\boldsymbol{P}\in\mathbb{R}^{d \times d_c}$, a relaxation scheme denoted as $\mathcal{G}$, residual tolerance $\delta$ and an initial approximation $\boldsymbol{u}^{(0)}$
\State set $k=0$, initial residual $\boldsymbol{r}^{(0)}=\boldsymbol{f}-\boldsymbol{K}\boldsymbol{u}^{(0)}$
\While{$\Vert \boldsymbol{r}^{(k)}\Vert < \delta$ }
    \State Pre-relaxation: Perform $r_1$ iterations of the relaxation scheme on the current approximation and obtain $\boldsymbol{u}^{(k)}$ as: $\boldsymbol{u}^{(k)} \gets \mathcal{G}\left(\boldsymbol{u}^{(k)};r_1 \right)$
    \State Update the residual: $\boldsymbol{r}^{(k)}=\boldsymbol{f}-\boldsymbol{K}\boldsymbol{u}^{(k)}$
    \State Restrict the residual to the coarser level and solve the coarse level system $\boldsymbol{K}_c\boldsymbol{e}_c^{(k)}=\boldsymbol{P}^T\boldsymbol{r}^{(k)}$, where $\boldsymbol{K}_c=\boldsymbol{P}^T \boldsymbol{K} \boldsymbol{P} \in \mathbb{R}^{d_c \times d_c}$
    \State Prolongate the coarse grid error $\boldsymbol{e}^{(k)}=\boldsymbol{P}\boldsymbol{e}_c^{(k)}$
    \State Correct the fine grid solution: $\boldsymbol{u}^{(k+1)}=\boldsymbol{u}^{(k)}+\boldsymbol{e}^{(k)}$
    \State Post-relaxation: Perform additional $r_2$ relaxation iterations and obtain $\boldsymbol{u}^{(k+1)} \gets \mathcal{G}\left(\boldsymbol{u}^{(k+1)};r_2 \right)$
    \State $k=k+1$
\EndWhile
\end{algorithmic}
\end{algorithm}

\begin{figure}[H]
\centering
    \begin{subfigure}{0.55\textwidth} 
      \includegraphics[width=\textwidth]{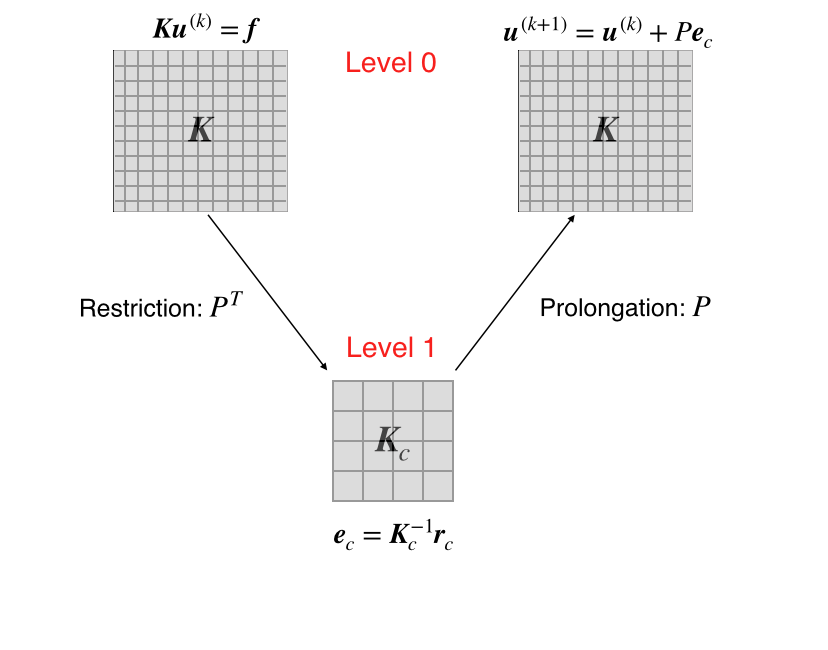}
      \subcaption{}
      \label{fig:2levelAMG}
     \end{subfigure}
    \begin{subfigure}{0.55\textwidth} 
      \includegraphics[width=\textwidth]{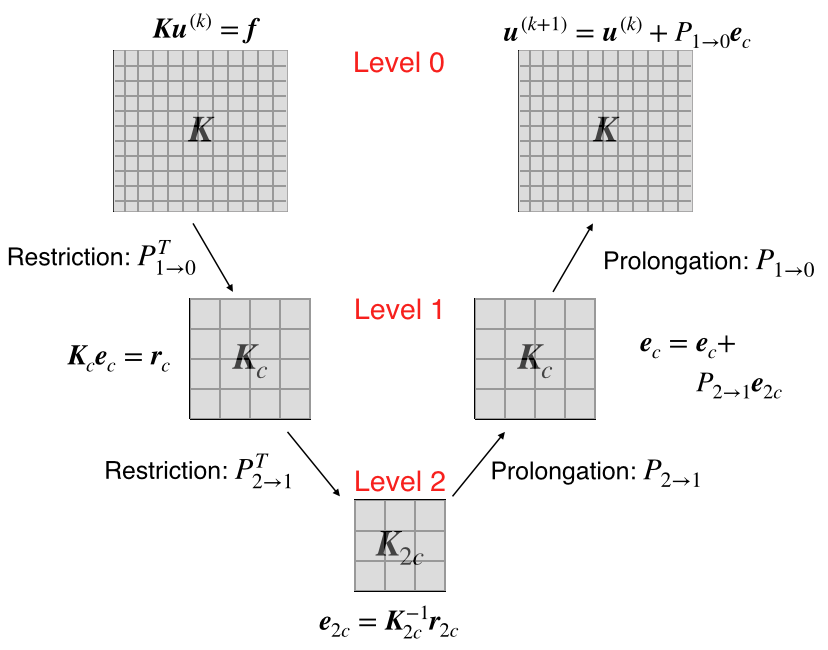}
      \subcaption{}
      \label{fig:3levelAMG}
    \end{subfigure}
\caption{Multigrid V-cycles in a (a) 2-level and a (b) 3-level setting}
\label{fig:1d_solutions}
\end{figure}

In the above algorithm, lines 4-10 describe what is known as a $V$-cycle, schematically depicted in figure \ref{fig:2levelAMG}. The multi-level version of the above algorithm is easily obtained as the result of recursively applying the two-level algorithm, as shown in fig. \ref{fig:3levelAMG} for the 3-level setting. The notation 

\begin{equation}
    \boldsymbol{u}^{(k+1)}=AMG(\boldsymbol{u}^{(k)};\boldsymbol{K},\boldsymbol{f},r_1,r_2)
\end{equation}

\noindent will be used to denote the application of one AMG cycle.

To better illustrate algorithm \ref{alg:AMGalgorithm} and its convergence properties, let us consider the GS algorithm as the relaxation scheme, where the matrix $\boldsymbol{K}$ is split into $\boldsymbol{K}=\boldsymbol{L}+\boldsymbol{V}$, $\boldsymbol{L}$ being a lower triangular matrix that includes the diagonal elements and $\boldsymbol{V}$ is the upper triangular part of $\boldsymbol{K}$. The iterative scheme of the GS method is as follows:

\begin{align} \label{eq:GaussSeidel}
    \boldsymbol{u}_{m+1}&=\boldsymbol{L}^{-1}\left(\boldsymbol{f}-\boldsymbol{V} \boldsymbol{u}_m\right) \nonumber \\
    &= \boldsymbol{L}^{-1}\boldsymbol{f}-\boldsymbol{L}^{-1}\left(\boldsymbol{K}-\boldsymbol{L} \right) \boldsymbol{u}_m \nonumber\\
    &= \boldsymbol{u}_{m}+\boldsymbol{L}^{-1}\left(\boldsymbol{f}-\boldsymbol{K}\boldsymbol{u}_m\right) \nonumber \\
    &= \boldsymbol{u}_{m}+\boldsymbol{L}^{-1}\boldsymbol{r}_m
\end{align}

\noindent where the subscripts $m,m+1$ in the above equation denote the iteration number of the GS algorithm. If $\boldsymbol{u}^\star$ is the exact solution to the system and $\boldsymbol{e}_m=\boldsymbol{u}^\star-\boldsymbol{u}_m$ the error after the $m$-th iteration, then

\begin{align} \label{eq:GaussSeidel_error}
    \boldsymbol{e}_{m+1}&=\boldsymbol{u}^{\star}-\boldsymbol{u}_{m+1} \nonumber \\
    &= \boldsymbol{e}_{m}+\boldsymbol{u}_{m}-\left(\boldsymbol{u}_{m}+\boldsymbol{L}^{-1}\boldsymbol{r}_m\right) \nonumber\\
    &= \boldsymbol{e}_{m}-\boldsymbol{L}^{-1}\left(\boldsymbol{K}\boldsymbol{e}_m\right) \nonumber\\
    &= \left(\boldsymbol{I}-\boldsymbol{L}^{-1}\boldsymbol{K} \right) \boldsymbol{e}_m
\end{align}

\noindent where $\boldsymbol{I}$ is the $d\times d$ identity matrix. Setting $\boldsymbol{M}=\boldsymbol{I}-\boldsymbol{L}^{-1}\boldsymbol{K}$, then it is straightforward to show that

\begin{equation}
    \boldsymbol{e}_{m+1}=\boldsymbol{M}\boldsymbol{e}_{m}=\boldsymbol{M}^2\boldsymbol{e}_{m-1}=\dots \boldsymbol{M}^{m+1}\boldsymbol{e}_{0}
\end{equation}

\noindent Returning to Algorithm \ref{alg:AMGalgorithm}, the error at the end of the $k$-th cycle of the two-level AMG can be computed as:

\begin{equation} \label{eq:AMGerror}
    \boldsymbol{e}^{(k)}= \boldsymbol{M}^{r_2}\boldsymbol{C}\boldsymbol{M}^{r_1}\boldsymbol{e}^{(k-1)}
\end{equation}

\noindent with

\begin{equation}
    \boldsymbol{C}= \boldsymbol{I}-\boldsymbol{P}\left( \boldsymbol{P}^T \boldsymbol{K} \boldsymbol{P} \right)^{-1}\boldsymbol{P}^T\boldsymbol{K} 
\end{equation}

\noindent being the coarse grid correction, $\boldsymbol{M}^{r_2}$ the post-relaxation matrix after $r_2$ sweeps and $\boldsymbol{M}^{r_1}$ the pre-relaxation after $r_1$ sweeps.

From eq. \eqref{eq:AMGerror} it becomes evident that the matrix $\boldsymbol{M}^{r_2}\boldsymbol{C}\boldsymbol{M}^{r_1}$ determines the convergence behavior of the two-level cycle. The relaxation matrix $\boldsymbol{M}$ plays a role, however, in practice the selection of the prolongation operator $\boldsymbol{P}$ is the key to designing an efficient algorithm. In this regard, the most popular variations of AMG include the Ruge-St{\"u}ben method \cite{Ruge1987} and the smoothed aggregation based (SA) AMG \cite{Vanek1996}. Lastly, another factor the affects the number of iterations in AMG to reach the prescribed threshold of accuracy, is the choice of the initial solution. In absence of other information, $\boldsymbol{u}^{(0)}=\boldsymbol{0}$ is usually considered.

\section{Machine learning accelerated iterative solvers} \label{sec3}

\subsection{Problem statement} \label{sec3.1}

The aim in this section is to develop an efficient data-driven and AI-enhanced solver for the parametrized system of eq. \eqref{eq:fineProblem}, by combining linear algebra-based solvers with machine learning algorithms. More specifically, the idea proposed herein, is to utilize a reduced set of high-fidelity system solutions, obtained after solving eq. \eqref{eq:fineProblem} for specified parameter instances, in two different yet complementary ways. First, a surrogate model will be established in the form of a `cheap-to-evaluate` nonlinear mapping from the problem's parameter space to its solution space using convolutional neural networks (CNNs) and feedforward neural networks (FFNNs). Even though CNNs and FFNNs have been shown to produce astonishing results even for challenging applications \cite{Mo2019, XU2020, NIKOLOPOULOS2022}, nevertheless, their black-box and physics-agnostic nature doesn't provide any means to improve the solutions they produce. To combat this problem, POD is performed on this data set of solutions and an efficient iterative solver is developed based on the idea of AMG, where in this case the prolongation operator is substituted by the projection matrix to the POD reduced space.

\subsection{Construction of surrogate model}\label{sec3.2} 
A surrogate model is an imitation of the original high fidelity model and serves as a 'cheap' mapping from the parametric space $\boldsymbol{\theta} \in \mathbb{R}^n$ to the solution space $\boldsymbol{u} \in \mathbb{R}^d$. In general, it is built upon an initial data set $\lbrace\boldsymbol{u}_i\rbrace_{i=1}^{N}$, which is created by solving the problem for a small, yet sufficient number, $N$, of parameter values. It is essential to span the problem's parametric space effectively, thus sophisticated sampling methods are often utilized, such as the Latin Hypercube \cite{olsson2002latin}. Many surrogate modeling techniques have been introduced over the past years, including linear \cite{Ladeveze2011,Zahr2017,Agathos2020} and nonlinear \cite{NIKOLOPOULOS2021, NIKOLOPOULOS2022, KALOGERIS2021} dimensionality reduction methods.  

In general, the selection of an appropriate surrogate modelling method is problem dependent, however, in this work, we will employ a surrogate modeling scheme based on convolutional autoencoders (CAEs) and feedforward neural networks (FFNNs) that was introduced in \cite{NIKOLOPOULOS2021} for parametrized time-dependent PDEs. It consists of two phases, namely the offline and the online phase. The offline phase begins with the training of a CAE that consists of an encoder and a decoder, in order to obtain low dimensional latent representations, $\boldsymbol{z}_{i} \in \mathbb{R}^l$ for each $\boldsymbol{u}_i \in \mathbb{R}^d$, through the encoder with $l\ll d$ and a reconstruction map by the decoder. It is trained over the initial data set $\lbrace\boldsymbol{u}_i\rbrace_{i=1}^{N}$ to minimize the objective function:

\begin{equation}\label{CAE}
\mathcal{L}_{CAE} = \frac{1}{N}\sum_{i=1}^{N}||\boldsymbol{u}_{i} - \tilde{\boldsymbol{u}}_{i}||_{2}^{2}
\end{equation}
where $\tilde{\boldsymbol{u}}_{i}$ is the reconstructed input. After the training is completed, the latent space data set $\lbrace\boldsymbol{z}_i\rbrace_{i=1}^{N}$ is obtained. The second step of the offline phase is the training of the FFNN, which is used to establish a nonlinear mapping from the parametric space $\boldsymbol{\theta} \in \mathbb{R}^n$  to the latent space $\boldsymbol{z} \in \mathbb{R}^l$. Again, the aim of the training is the minimization of the loss function:

\begin{equation}\label{FFNN}
\mathcal{L}_{FFNN} = \frac{1}{N}\sum_{i=1}^{N}||\boldsymbol{z}_{i} - \tilde{\boldsymbol{z}}_{i}||_{2}^{2}
\end{equation}
where $\tilde{\boldsymbol{z}}_{i}$ is the network's output.

Subsequently, the online phase utilizes the fully trained surrogate model, which is now capable of delivering accurate predictions of the system's response for new parameter values $\boldsymbol{\theta}_{j}$ as follows:

\begin{equation}
\boldsymbol{u}_{j} = decoder(FFNN(\boldsymbol{\theta}_{j})):=\mathcal{F}^{sur}(\boldsymbol{\theta}_j)
\end{equation}

\begin{figure}[H]
    \centering
    \includegraphics[width=0.7\textwidth]{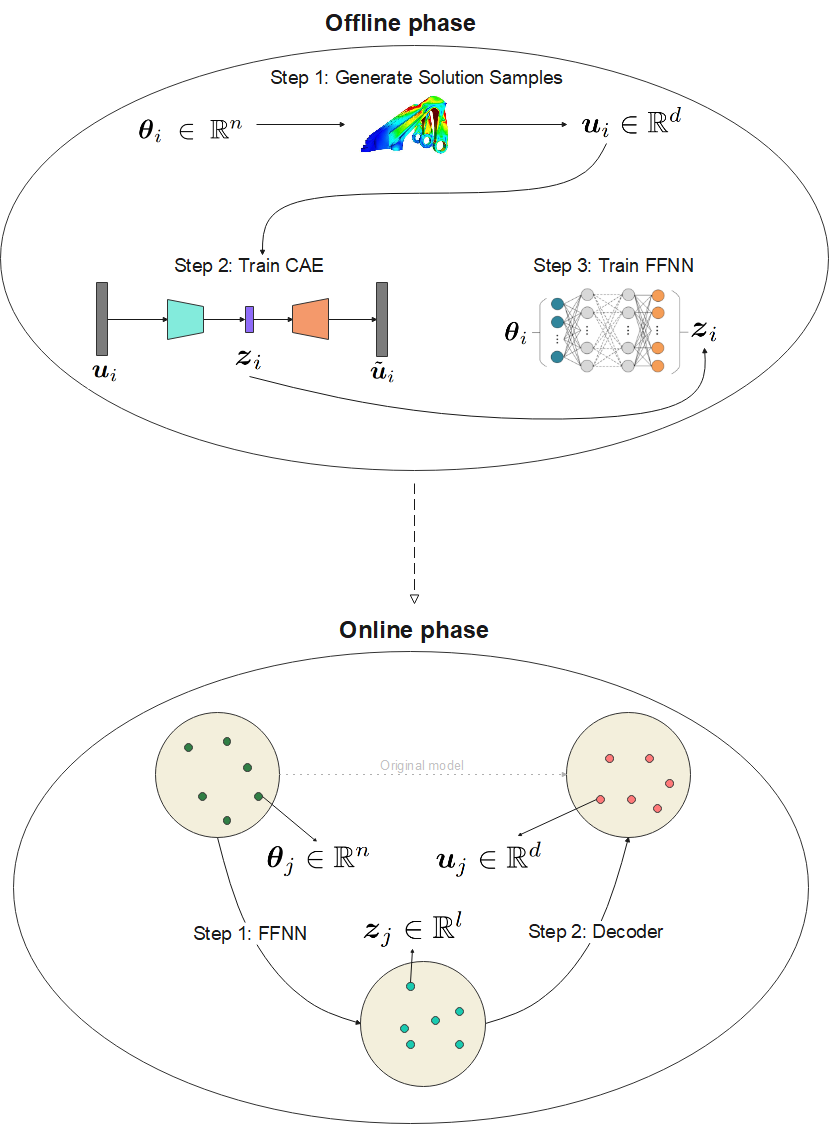}
    \caption{Schematic representation of the surrogate model}
    \label{fig:Surrogate}
\end{figure}

A schematic representation of the surrogate model is presented in figure \ref{fig:Surrogate}.

\subsection{Multigrid-inspired POD solver}\label{sec3.3} 
POD, also known as Principal Component Analysis, is a powerful and effective approach for data analysis and dimensionality reduction, aimed at indentifying low-order modes of a system. In conjunction with the Galerkin projection procedure it is commonly utilized as an efficient method to reduce the dimensionality of large linear systems of equations \cite{Rathinam2003, LIEU2006, RAPUN2010}. The theory and application of POD is covered in many publications, however, to keep this paper as self-contained as possible the POD procedure used within this framework is summarized below. Let us denote with $\boldsymbol{U}\in \mathbb{R}^{d\times N}$ the matrix consisting of $N$ solution vectors $\left[\boldsymbol{u}_1,...,\boldsymbol{u}_N \right]$ for different parameter values $\lbrace \boldsymbol{\theta}_i\rbrace_{i=1}^N$ and with $\boldsymbol{R}=\boldsymbol{U}\boldsymbol{U}^T\in \mathbb{R}^{d \times d}$ the correlation matrix. Then POD consists in the following steps.

\begin{enumerate}
    \item Compute the eigenvalues and eigenvectors of $\boldsymbol{R}$ that satisfy $\boldsymbol{R} \boldsymbol{\Phi}=\boldsymbol{\Phi}\boldsymbol{\Lambda}$. This step can be very demanding when $d\gg 1$, however, in practice $N\ll d$ and since $\boldsymbol{R}$, $\boldsymbol{R}^T$ have the same non-zero eigenvalues, it is computationally more convenient to solve instead the eigenvalue problem $\boldsymbol{U}^T\boldsymbol{U} \boldsymbol{\Psi}=\boldsymbol{\Psi}\boldsymbol{\Lambda}$. Then, the eigenvectors $\boldsymbol{\Phi}$ and $\boldsymbol{\Psi}$ are linked according to the formula.

    \begin{equation} \label{eq:PODrelationEigenvectors}
    \boldsymbol{\Phi}=\boldsymbol{U}\boldsymbol{\Psi}\boldsymbol{\Lambda}^{-1/2}
    \end{equation}

    \item Form the reduced basis $\boldsymbol{\Phi}_r$, be retaining only the $r$ first columns of $\boldsymbol{\Phi}$, corresponding to the largest eigenvalues.
    
    \item Under the assumption that each solution to eq. \eqref{eq:fineProblem} can be approximated as:
    
    \begin{equation}
        \boldsymbol{u}\equiv \boldsymbol{\Phi}_r \boldsymbol{u}_r
    \end{equation}
    
    \noindent with $\boldsymbol{u}_r\in \mathbb{R}^r$ being the unknown coefficients of the projection on the truncated POD basis, then the reduced-order linear system becomes:
    
    \begin{align} \label{eq:POD_ROM}
        &\boldsymbol{K}\boldsymbol{u}=\boldsymbol{f} \nonumber \\
        &\boldsymbol{\Phi}_r^T\boldsymbol{K}\boldsymbol{\Phi}_r \boldsymbol{u}_r=\boldsymbol{\Phi}_r^T \boldsymbol{f} \nonumber \\
        &\boldsymbol{K}_r \boldsymbol{u}_r=\boldsymbol{f}_r
    \end{align}
    
    \noindent Solving equation \eqref{eq:POD_ROM} for $\boldsymbol{u}_r$ is significantly easier since $\boldsymbol{K}_r \in \mathbb{R}^{r \times r}$, with $r$ small. 
    
    \item Retrieve the solution to their original problem:
    
    \begin{equation}
        \boldsymbol{u}=\boldsymbol{\Phi}_r\boldsymbol{u}_r
    \end{equation}
    
\end{enumerate}

Based on the above, a similarity between the 2-level AMG method and POD can be observed, under the identification of $\boldsymbol{\Phi}_r$ as the prolongation operator and $\boldsymbol{\Phi}_r^T$ the corresponding restriction. Then, the algorithm \ref{alg:AMGalgorithm} remains practically the same. In this case, the error of the scheme is given by the formula

\begin{equation} \label{eq:FinalScheme}
    \boldsymbol{e}^{(k)}=\boldsymbol{M}^{r_2}\left( \boldsymbol{I}- \boldsymbol{\Phi}_r\left( \boldsymbol{\Phi}_r^T\boldsymbol{K}\boldsymbol{\Phi}_r\right)^{-1} \boldsymbol{\Phi}_r^T\boldsymbol{K} \right)\boldsymbol{M}^{r_1}\boldsymbol{e}^{(k-1)}
\end{equation}

\subsection{Proposed data-driven framework for parameterized linear systems} \label{sec3.4}

The final step is to combine the surrogate model of section \ref{sec3.2} and the multigrid-inspired POD solver of the previous section into a unified methodological framework for solving efficiently large-scale parametrized linear systems. In particular, an initial data set of system solutions $\lbrace\boldsymbol{u}_i\rbrace_{i=1}^{N}$ is performed for specified instances of the parameter vector $\lbrace\boldsymbol{\theta}_i\rbrace_{i=1}^{N}$. Then, these solution vectors are utilized as training data for the CNN and FFNN and the surrogate model is established. The error between the exact solution and the surrogate's prediction for a given $\boldsymbol{\theta}$ can be given as:

\begin{equation}
    \boldsymbol{e}^{sur}=\boldsymbol{u}^{\star} -\mathcal{F}^{sur}(\boldsymbol{\theta})
\end{equation}

Despite one's best efforts, however, $\Vert \boldsymbol{e}^{sur}\Vert \neq 0$ and the surrogate's predictions will not satisfy exactly equation \eqref{eq:fineProblem}. At this point, instead of simply performing iterations of PCG or AMG to improve the surrogate's predictions, we propose to further utilize the knowledge available to us from the data set of solution vectors, in order to enhance the performance of these iterative solvers. In particular, we perform POD to the solution matrix $\boldsymbol{U}=[\boldsymbol{u}_1,...,\boldsymbol{u}_{N}]$, in order to obtain the projection matrix $\boldsymbol{\Phi}_r^T$ and apply the AMG method either directly, or as a preconditioner in the PCG algorithm according to the following algorithm \ref{alg:PCG_AMGprecondAlgorithm}.  

\begin{algorithm} [H]
\setstretch{1.0}
\caption{AMG preconditioned PCG algorithm}\label{alg:PCG_AMGprecondAlgorithm}
\begin{algorithmic}[1]
\State \textbf{Input:} $\boldsymbol{K}\in\mathbb{R}^{d\times d}$, rhs $\boldsymbol{f}\in\mathbb{R}^d$, AMG scheme, residual tolerance $\delta$ and an initial approximation $\boldsymbol{u}^{(0)}$
\State set $k=0$, initial residual $\boldsymbol{r}^{(0)}=\boldsymbol{f}-\boldsymbol{K}\boldsymbol{u}^{(0)}$
\State $\boldsymbol{s}_0=AMG(\boldsymbol{0};\boldsymbol{K},\boldsymbol{r}^{(0)},r_1,r_2)$
\State $\boldsymbol{p}_0=\boldsymbol{s}_0$
\While{$\Vert \boldsymbol{r}^{(k)}\Vert < \delta$ }
    \State $\alpha_k=\frac{\left(\boldsymbol{r}^{(k)}\right)^T\boldsymbol{s}_k}{\boldsymbol{p}_k^T\boldsymbol{K}\boldsymbol{p}_k}$
    \State $\boldsymbol{u}^{(k+1)}=\boldsymbol{u}^{(k)}+\alpha_k \boldsymbol{p}_k$
    \State $\boldsymbol{r}^{(k+1)}=\boldsymbol{r}^{(k)}-\alpha_k \boldsymbol{K}\boldsymbol{p}_k$
    \State $\boldsymbol{s}_{k+1}=AMG(\boldsymbol{0};\boldsymbol{K},\boldsymbol{r}^{(k+1)},r_1,r_2)$
        \State $\beta_k=\frac{\left(\boldsymbol{r}^{(k+1)}\right)^T\boldsymbol{s}_{k+1}}{\left(\boldsymbol{r}^{(k)}\right)^T\boldsymbol{s}_{k}}$
    \State $\boldsymbol{p}_{k+1}=\boldsymbol{s}_{k+1}+\beta_k\boldsymbol{p}_k$
    \State $k=k+1$
\EndWhile
\end{algorithmic}
\end{algorithm}

\subsection{Error bounds} \label{sec3.5}

The proposed methodology is data-driven and, as such, it is not possible to provide a priori estimates of the error for general systems. Nevertheless, under the assumption that the training data set $\boldsymbol{U}$ is `large' enough to contain almost all possible variations of the solution vector, then an estimate for the error can be provided as follows:

\begin{align} \label{ineq:Main}
    \boldsymbol{e}^{(k)}&=\boldsymbol{M}^{r_2}\left( \boldsymbol{I}- \boldsymbol{\Phi}_r\left( \boldsymbol{\Phi}_r^T\boldsymbol{K}\boldsymbol{\Phi}_r\right)^{-1} \boldsymbol{\Phi}_r^T\boldsymbol{K} \right)\boldsymbol{M}^{r_2}\boldsymbol{e}^{(k-1)} \Rightarrow \nonumber \\
    \Vert \boldsymbol{e}^{(k)} \Vert &=\Vert \boldsymbol{M}^{r_2}\left( \boldsymbol{I}- \boldsymbol{\Phi}_r\left( \boldsymbol{\Phi}_r^T\boldsymbol{K}\boldsymbol{\Phi}_r\right)^{-1} \boldsymbol{\Phi}_r^T\boldsymbol{K} \right)\boldsymbol{M}^{r_2}\boldsymbol{e}^{(k-1)} \Vert \Rightarrow \nonumber \\
    &\leq \Vert \boldsymbol{M}^{r_2} \Vert \Bigg{\Vert} \left( \boldsymbol{I}- \boldsymbol{\Phi}_r\left( \boldsymbol{\Phi}_r^T\boldsymbol{K}\boldsymbol{\Phi}_r\right)^{-1} \boldsymbol{\Phi}_r^T\boldsymbol{K} \right) \Bigg{\Vert} \Vert \boldsymbol{M}^{r_2}\Vert \Vert \boldsymbol{e}^{(k-1)} \Vert
\end{align}

\noindent In the above, $\Vert \cdot \Vert$ denotes the $l2$-vector norm, when the input is a vector, and the induced operator norm (spectral norm) when the input is a matrix.

We can assume that $\boldsymbol{M}$ has a spectral radius $\rho(\boldsymbol{M})<1$ and the GS algorithm converges. This assumption is valid when $\boldsymbol{K}$ is symmetric positive definite, which is commonly the case in engineering problems. Then, according to Gelfand's formula, we have

\begin{equation}
    \rho(\boldsymbol{M})=\lim_{k\rightarrow \infty}\Vert \boldsymbol{M}^k \Vert^{1/k}
\end{equation}

\noindent As a consequence, there is $k_0\in\mathbb{N}$ and $\gamma\in \left(\rho\left(\boldsymbol{M}\right),1\right)\subseteq (0,1)$
such that:

\begin{equation} \label{ineq:sideTerms}
    \Vert \boldsymbol{M}^k \Vert \leq \gamma^k, \ \ \ \forall k \geq k_0
\end{equation}

\noindent Therefore,  $\Vert \boldsymbol{M}^{r_1} \Vert,\Vert \boldsymbol{M}^{r_2} \Vert < 1$ for $r_1, r_2$ large enough. 

Now, focusing on the term $\Bigg{\Vert} \left( \boldsymbol{I}- \boldsymbol{\Phi}_r\left( \boldsymbol{\Phi}_r^T\boldsymbol{K}\boldsymbol{\Phi}_r\right)^{-1} \boldsymbol{\Phi}_r^T\boldsymbol{K} \right) \Bigg{\Vert}$, then, by definition the following holds:

\begin{equation}
    \Bigg{\Vert} \left( \boldsymbol{I}- \boldsymbol{\Phi}_r\left( \boldsymbol{\Phi}_r^T\boldsymbol{K}\boldsymbol{\Phi}_r\right)^{-1} \boldsymbol{\Phi}_r^T\boldsymbol{K} \right) \Bigg{\Vert}=\sup_{\boldsymbol{u}\in\mathbb{R}^d : \Vert \boldsymbol{u}\Vert = 1}  \Bigg{\Vert} \left( \boldsymbol{I}- \boldsymbol{\Phi}_r\left( \boldsymbol{\Phi}_r^T\boldsymbol{K}\boldsymbol{\Phi}_r\right)^{-1} \boldsymbol{\Phi}_r^T\boldsymbol{K} \right)\boldsymbol{u} \Bigg{\Vert}
\end{equation}

\noindent Since a given $\boldsymbol{u} \in \mathbb{R}^d$ can be decomposed as $\boldsymbol{u}=\boldsymbol{\Phi}_r\boldsymbol{u}_r + \boldsymbol{u}^\perp$, with $\boldsymbol{\Phi}_r\boldsymbol{u}_r\in \Phi$ and  $\boldsymbol{u}^\perp\in \Phi^\perp$, where $\Phi=\text{span } \lbrace \boldsymbol{\phi}_1,...,\boldsymbol{\phi}_r \rbrace $ and $\Phi^\perp$ its orthogonal complement in $\mathbb{R}^d$, then,

\begin{align}
    \left( \boldsymbol{I}- \boldsymbol{\Phi}_r\left( \boldsymbol{\Phi}_r^T\boldsymbol{K}\boldsymbol{\Phi}_r\right)^{-1} \boldsymbol{\Phi}_r^T\boldsymbol{K} \right)\boldsymbol{u} &= \boldsymbol{u}- \boldsymbol{\Phi}_r\left( \boldsymbol{\Phi}_r^T\boldsymbol{K}\boldsymbol{\Phi}_r\right)^{-1} \boldsymbol{\Phi}_r^T\boldsymbol{K} \left(\boldsymbol{\Phi}_r\boldsymbol{u}_r+\boldsymbol{u}^\perp\right) \nonumber \\
    &= \boldsymbol{\Phi}_r\boldsymbol{u}_r+\boldsymbol{u}^\perp- \boldsymbol{\Phi}_r\left( \boldsymbol{\Phi}_r^T\boldsymbol{K}\boldsymbol{\Phi}_r\right)^{-1} \boldsymbol{\Phi}_r^T\boldsymbol{K} \left(\boldsymbol{\Phi}_r\boldsymbol{u}_r+\boldsymbol{u}^\perp\right) \nonumber \\
    &= \boldsymbol{\Phi}_r\boldsymbol{u}_r+\boldsymbol{u}^\perp- \boldsymbol{\Phi}_r\boldsymbol{u}_r-\boldsymbol{\Phi}_r\left( \boldsymbol{\Phi}_r^T\boldsymbol{K}\boldsymbol{\Phi}_r\right)^{-1} \boldsymbol{\Phi}_r^T\boldsymbol{K} \boldsymbol{u}^\perp \nonumber \\
    &= \boldsymbol{u}^\perp-\boldsymbol{\Phi}_r\left( \boldsymbol{\Phi}_r^T\boldsymbol{K}\boldsymbol{\Phi}_r\right)^{-1} \boldsymbol{\Phi}_r^T\boldsymbol{K} \boldsymbol{u}^\perp 
\end{align}

\noindent thus,

\begin{equation} \label{ineq:middleTerm}
        \Bigg{\Vert} \left( \boldsymbol{I}- \boldsymbol{\Phi}_r\left( \boldsymbol{\Phi}_r^T\boldsymbol{K}\boldsymbol{\Phi}_r\right)^{-1} \boldsymbol{\Phi}_r^T\boldsymbol{K} \right)\boldsymbol{u}\Bigg{\Vert}  \leq  \Bigg{\Vert} \boldsymbol{I}-\boldsymbol{\Phi}_r\left( \boldsymbol{\Phi}_r^T\boldsymbol{K}\boldsymbol{\Phi}_r\right)^{-1} \boldsymbol{\Phi}_r^T\boldsymbol{K}\Bigg{\Vert} \Vert \boldsymbol{u}^\perp \Vert \leq c\Vert \boldsymbol{u}^\perp \Vert 
\end{equation}

\noindent for some $c > 0$. Due to the orthogonality of $\boldsymbol{\Phi}_r\boldsymbol{u}_r$ and $\boldsymbol{u}^\perp $, it follows that 

\begin{align}
    &\Vert \boldsymbol{u}^\perp \Vert ^2 = \Vert \boldsymbol{u} \Vert ^2 - \Vert \boldsymbol{\Phi}_r\boldsymbol{u}_r  \Vert ^2 \Rightarrow \nonumber \\
    &\Vert \boldsymbol{u}^\perp \Vert = \sqrt{ 1 - \Vert \boldsymbol{\Phi}_r\boldsymbol{u}_r  \Vert ^2} \leq 1
\end{align}

\noindent In fact, by choosing an appropriate number of eigenvectors $r$ in POD, we can obtain $\Vert \boldsymbol{u}^\perp \Vert < \frac{1}{c}$ and then, inequality \eqref{ineq:middleTerm} becomes

\begin{equation} \label{ineq:middleTermUpdate}
        \Bigg{\Vert} \left( \boldsymbol{I}- \boldsymbol{\Phi}_r\left( \boldsymbol{\Phi}_r^T\boldsymbol{K}\boldsymbol{\Phi}_r\right)^{-1} \boldsymbol{\Phi}_r^T\boldsymbol{K} \right)\boldsymbol{u}\Bigg{\Vert}  \leq  C
\end{equation}

\noindent with $C\coloneqq C(\boldsymbol{u}^\perp)$ and $C\in (0,1)$.

Inserting the inequalities \eqref{ineq:sideTerms} and \eqref{ineq:middleTermUpdate} into \eqref{ineq:Main}, we have:

\begin{equation}
    \Vert \boldsymbol{e}^{(k)} \Vert \leq \gamma^{r_2} C \gamma^{r_1} \Vert \boldsymbol{e}^{(k-1)} \Vert, \quad \text{with  }  \gamma^{r_2} C \gamma^{r_1}<1 
\end{equation}

\noindent Applying the above inequality recursively, we conclude:

\begin{align}
    \Vert \boldsymbol{e}^{(k)} \Vert &\leq (\gamma^{r_2})^k C^k  (\gamma^{r_1})^k \Vert \boldsymbol{e}^{(0)} \Vert   \nonumber \\
    &=(\gamma^{r_2})^k C^k  (\gamma^{r_1})^k \Vert \boldsymbol{e}^{sur} \Vert
\end{align}

The above inequality provides us with some valuable insight regarding the performance of the proposed data-driven solver. Most importantly, we notice the critical role that the surrogate's predictions play in the convergence, since the error is bounded be the surrogate's error $\Vert \boldsymbol{e}^{sur} \Vert$. Even though this result agrees with common intuition, nevertheless, being rigorously proven excludes the possibility of good initial predictions requiring more iterations for the solution to converge. Secondly, by retaining more eigenvectors to construct the reduced space $\Phi$, we reduce the norm of $\boldsymbol{u}^\perp\in \Phi^\perp $, resulting in faster convergence. In the following section, we test the solver on numerical applications of scientific interest and assess its performance in comparison with conventional solvers.

\section{Numerical applications} \label{sec4}
The proposed methodology is tested on two large scale parametrized systems. The first case is the indirect tensile strength (ITS) test, which is treated with the theory of 2D linear elasticity, while the second one is a 3D deformable porous medium problem, also known as Biot problem.

\subsection{Indirect tensile strength test} 

\begin{figure}[H]
    \centering
    \includegraphics[width=0.5\textwidth]{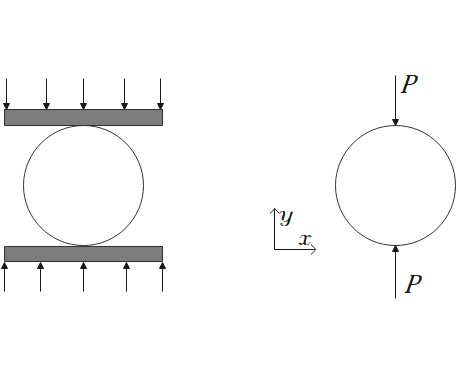}
    \caption{ITS test: A diametrically point loaded disk}
    \label{fig:ex1}
\end{figure}

A popular test to measure the tensile strength of concrete or asphalt materials is the ITS test. As shown in figure \ref{fig:ex1}, the test contains a cylindrical specimen loaded across its diameter to failure. The specimen is usually loaded at a constant deformation rate and measuring the load response. When the developed tensile stress in the specimen under loading exceeds its tensile strength then the specimen will fail. In this application, we restrict our analysis to the linear regime and model the cylinder as a 2D disk under plain strain assumptions, as shown in \ref{fig:ex1}. In this case, the weak form of the problem reads: Find $\boldsymbol{v}\in\mathcal{V}(\Omega)$ such that

\begin{equation}
\begin{aligned}
     \int_{\Omega}& \boldsymbol{\sigma}(\boldsymbol{v}): \boldsymbol{\epsilon}(\boldsymbol{w})d\Omega
    = \int_{\Omega}\boldsymbol{f}\cdot \boldsymbol{w}d\Omega, \quad  \forall \boldsymbol{w}\in\mathcal{V}_c(\Omega)  \\
   \boldsymbol{\sigma} &= \lambda tr\left(\boldsymbol{\epsilon} \right)\mathbb{I}+2\mu\boldsymbol{\epsilon}
\end{aligned} 
\label{eq:ITS}
\end{equation}
where, 

\begin{equation}
    \boldsymbol{\epsilon}=\begin{bmatrix} \epsilon_{xx} & \epsilon_{xy} \\  \epsilon_{xy} & \epsilon_{yy} 
    \end{bmatrix}
\end{equation}

\noindent the strain tensor and $\boldsymbol{f}$ the loading. Also, $\mu$ and $\lambda$ are the Lamé's constants, which are linked to the Young modulus $E$ and the Poisson ratio according to equations \eqref{eq:Lame}:

\begin{equation}
\begin{aligned}
\mu &= \frac{E}{2(1 + \nu)} \\
\lambda &= \frac{E\nu}{(1 + \nu)(1 - 2\nu)}
\end{aligned} 
\label{eq:Lame}
\end{equation}

\begin{figure}[H]
    \centering
    \includegraphics[width=0.65\textwidth]{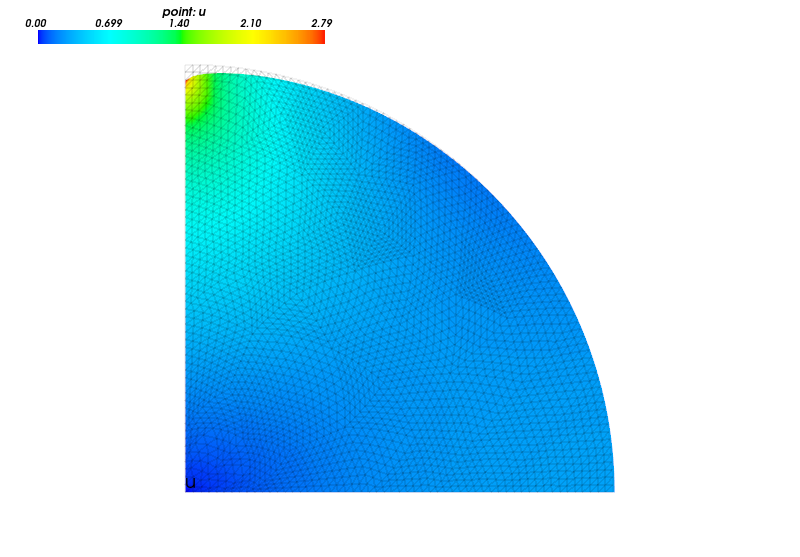}
    \caption{Displacement magnitude $||\boldsymbol{u}||$ for $E = 2000 \ MPa$ and $P = -1000 \ N$}
    \label{fig:snapshot}
\end{figure}

In this example, the specimen has a diameter of $150 \ mm$ and due to symmetry in geometry and loading we only need to model one quarter of the disk, as illustrated in figure \ref{fig:snapshot}. The solution of eq. \eqref{eq:ITS} is obtained using a finite element mesh that consists of triangular plane-strain finite elements with a total of $d = 5656$ dofs. The Young modulus $E$ and the load  $P$ are considered uncorrelated random variables following the lognormal distribution as described in table \ref{table:RandomParameters}. The Poisson ratio is considered to be a constant parameter $\nu = 0.3$. Figure \ref{fig:snapshot} displays the contour plot of the displacement norm $ \Vert \boldsymbol{u} \Vert$ for the mean value of the random parameters, that is $E = 2000 \ MPa$ and $P = -1000 \ N$.

\begin{table}[H]
 \centering
\begin{tabular}{|c|c|c|c|} 
\hline
 Parameter & Distribution & Mean & Standard deviation \\
\hline
$E(MPa)$ & Lognornal & $2000$ & $600$\\ 
\hline
$P(N)$ & Lognormal & $-1000$ & $300$\\ 
\hline
\end{tabular}
\caption{Random parameters of the ITS test}
\label{table:RandomParameters}
\end{table}

\begin{figure}[H]
    \centering
    \includegraphics[width=1\textwidth]{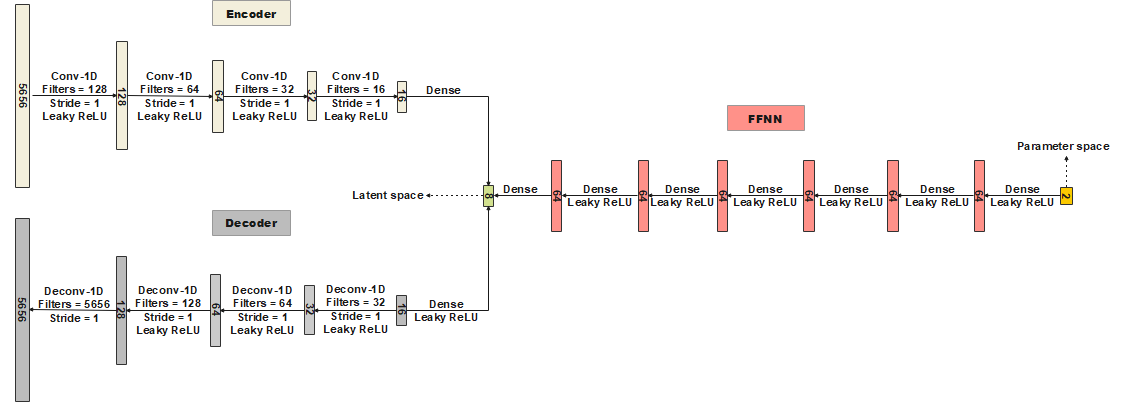}
    \caption{Surrogate model architecture}
    \label{fig:surrogate_1}
\end{figure}

The first step of the proposed procedure is to generate a sufficient number of offline samples. To this purpose, the Latin Hybercube sampling method was utilized to generate $N = 200$ parameter samples $\{[E_i, P_i]\}_{i=1}^{N}$. Subsequently, the corresponding problems are solved with the finite element method and the solution vectors obtained, $\{\boldsymbol{u}_i\}_{i=1}^{N}$, are regarded as 'exact' solutions. 
Next, a surrogate model is trained over these solutions in order to establish a 'cheap' mapping from the parametric to solution space.  The methodology for the surrogate model is described in section \ref{sec3.1}. The details of the selected CAE and FFNN architectures are presented in figure \ref{fig:surrogate_1}. 

To tackle the problem of overfitting, the standard hold-out approach was employed. In particular, the data set was randomly divided into train and validation subsets with a ratio of 70\%-30\% and each network's performance on the validation data set was assessed in order to avoid overfitting. The CAE is trained for 40 epochs with a batch size of 10 and a learning rate of $0.0005$, while the FFNN is trained for 3000 epochs with a batch size of 20 and a learning rate of $0.0001$. The average normalized $l_2$ norm error of the surrogate model on the test data set is 0.54\%.

The second step is to form the POD basis $\boldsymbol{\Phi}_r$ by performing eigendecomposition on the correlation matrix $\boldsymbol{U}\boldsymbol{U}^T$, with $\boldsymbol{U} = [\boldsymbol{u}_1,\dots,\boldsymbol{u}_N]$ being the solution matrix. In this case, the number of eigenvectors kept is $r=8$, which correspond to over 99.99\% of the variance in the training data. Subsequently, when all components of the proposed POD-2G solver are defined and fully trained, the methodology described in section \ref{sec3.2} can be applied to obtain new system's solutions for different parameter values. 

\begin{figure}[H]
\centering
\begin{subfigure}[b]{0.49\textwidth}
         \centering
         \includegraphics[width=\textwidth]{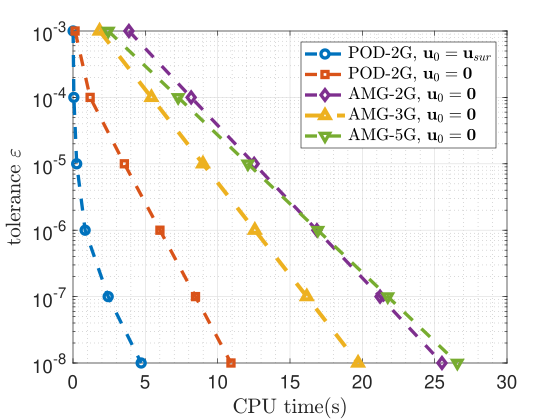}
         \caption{Comparison of mean CPU time}
     \end{subfigure}
     \hfill
     \begin{subfigure}[b]{0.49\textwidth} \label{fig:Results_AMG_1_subfigb}
         \centering
         \includegraphics[width=\textwidth]{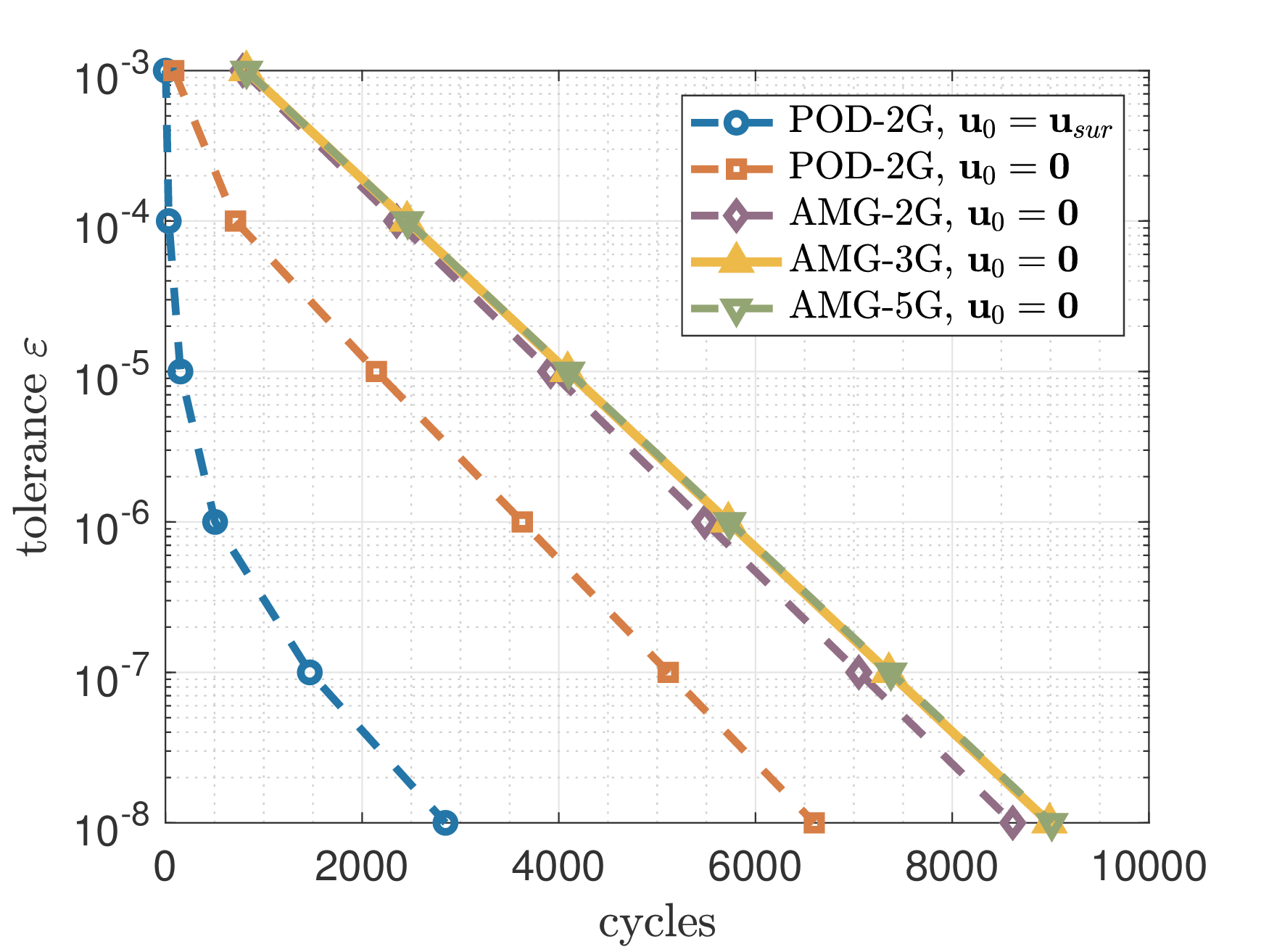}
          \caption{Comparison of mean number of cycles}
     \end{subfigure} 
         \caption{Comparison of mean CPU time and mean number of cycles over 500 analyses for different multigrid solvers}
    \label{fig:Results_AMG_1}
 \end{figure}

In order to test the proposed POD-based solver, a number of $N_{test} = 500$ test parameter samples $\{[E_j, P_j]\}_{j=1}^{N_{test}}$ were generated according to their distribution. The corresponding problems were solved with the Ruge-St{\"u}ben AMG solver for 2, 3 and 5 grids (termed AMG-2G, -3G, -5G respec.), as well as the proposed POD-2G solver for different values of tolerance. The size of the system of equations at the coarsest level for each of these solvers is given in table \ref{table:ProblemSizeEx1}. The mean value of the CPU time and the number of cycles required for convergence to the desired tolerance are displayed in figure \ref{fig:Results_AMG_1} for the 3 AMG solvers and the proposed POD-2G with initial $\boldsymbol{u}^{(0)}=\boldsymbol{0}$, as well as $\boldsymbol{u}^{(0)}=\boldsymbol{u}_{sur}$, namely the solution delivered by the surrogate model. 

\begin{table}[H]
 \centering
\begin{tabular}{|c|c|} 
\hline
 {} & System Size  \\
\hline
Initial Problem & $5656\times 5656$\\ 
\hline
AMG-2G & $1555 \times 1555$\\ 
\hline
AMG-3G & $314 \times 314$\\ 
\hline
AMG-4G & $80 \times 80$\\ 
\hline
AMG-5G & $26 \times 26$\\ 
\hline
POD-2G & $8 \times 8$\\ 
\hline
\end{tabular}
\caption{Size of the problem at the coarsest grid for the different solvers}
\label{table:ProblemSizeEx1}
\end{table}

\begin{table}[H]
 \centering
\begin{tabular}{|c||c|c|c|c|c|} 
\hline
  & $\varepsilon=10^{-4}$ & $\varepsilon=10^{-5}$ & $\varepsilon=10^{-6}$ & $\varepsilon=10^{-7}$ & $\varepsilon=10^{-8}$ \\
\hhline{|=||=|=|=|=|=|}
AMG-2G ($\boldsymbol{u}^{(0)}=\boldsymbol{0}$) & $\times$1 & $\times$1 & $\times$1 & $\times$1 & $\times$1\\ 
\hline
AMG-3G ($\boldsymbol{u}^{(0)}=\boldsymbol{0}$) & $\times$1.51 & $\times$1.39 & $\times$1.34 & $\times$1.31 & $\times$1.29\\ 
\hline
AMG-5G ($\boldsymbol{u}^{(0)}=\boldsymbol{0}$) & $\times$1.12 & $\times$1.04 & $\times$1.00 & $\times$0.98 & $\times$0.96\\ 
\hline
POD-2G ($\boldsymbol{u}^{(0)}=\boldsymbol{0}$) & $\times$6.90 & $\times$3.53 & $\times$2.81 & $\times$2.51 & $\times$2.34\\
\hline
POD-2G ($\boldsymbol{u}^{(0)}=\boldsymbol{u}_{sur}$) & $\times$138.99 & $\times$48.97 & $\times$20.09 & $\times$8.73 & $\times$5.51\\
\hline
\end{tabular}
\caption{Computational speedup of different solvers compared to AMG-2G}
\label{table:EX1_ComparisonAMG}
\end{table}
 
From figure \ref{fig:Results_AMG_1}, we notice that AMG-3G and AMG-5G require the same mean number of cycles, which is slightly more than those AMG-2G needs to achieve the same tolerance, yet, AMG-3G is the most efficient AMG scheme in terms of CPU time. This is because the CPU time is affected by both the size of the coarse problem and the number of times the prolongation/restriction operators are applied within a cycle. In this regard, AMG-3G provides the optimal number of grids needed for this problem. However, a significant improvement on both the speedup and the number of iterations can be observed when applying the two POD solvers instead of the AMG solvers (see table \ref{table:EX1_ComparisonAMG}). This performance gain is increased with increasing tolerance $\varepsilon$, reaching a speedup of $\times 6.90$ and $\times 138.99$ for the POD solvers with $\boldsymbol{u}^{(0)}=\boldsymbol{0}$ and $\boldsymbol{u}^{(0)}=\boldsymbol{u}_{sur}$ for $\varepsilon=10^{-4}$, respectively. On the other hand, for smaller values of $\varepsilon$ such as $10^{-8}$, the speedup in CPU time obtained with POD-2G with $\boldsymbol{u}^{(0)}=\boldsymbol{u}_{sur}$ is $\times 5.51$, when compared with the case of AMG-2G with $\boldsymbol{u}^{(0)} = \boldsymbol{0}$. Even though the gain achieved in this case is much smaller than for the case of $\varepsilon=10^{-4}$, yet, it is still notable. Based on these results, the conclusion is drawn that a key component of the proposed methodology is to obtain a close estimation of the solution by the surrogate model, $\boldsymbol{u}^{(0)} = \boldsymbol{u}_{sur}$ since an initial solution $\boldsymbol{u}^{(0)}$ from an accurately trained surrogate is capable of drastically reducing the computational cost.

\begin{figure}[H]
\centering
\begin{subfigure}[b]{0.49\textwidth}
         \centering
         \includegraphics[width=\textwidth]{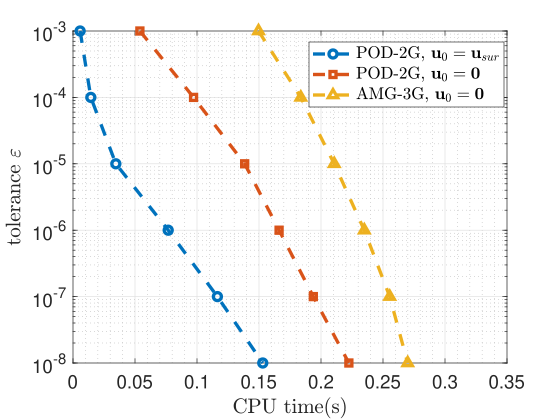}
         \caption{Comparison of mean CPU time}
     \end{subfigure}
     \hfill
     \begin{subfigure}[b]{0.49\textwidth}
         \centering
         \includegraphics[width=\textwidth]{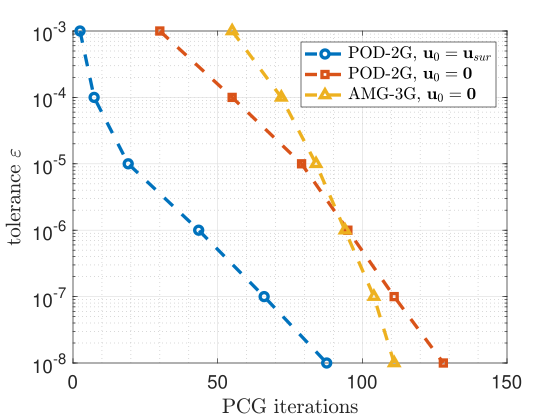}
          \caption{Comparison of mean PCG iterations}
     \end{subfigure} 
         \caption{Comparison of mean CPU time and mean number of PCG iterations over 500 analyses for different preconditioners}
    \label{fig:Results_PCG_1}
 \end{figure}

Furthermore, the convergence behaviour of the proposed method when used as a preconditioner in the context of the PCG method is presented in figure \ref{fig:Results_PCG_1} and table \ref{table:EX1_ComparisonPCG}.

 \begin{table}[H]
 \centering
\begin{tabular}{|c||c|c|c|c|c|} 
\hline
  & $\varepsilon=10^{-4}$ & $\varepsilon=10^{-5}$ & $\varepsilon=10^{-6}$ & $\varepsilon=10^{-7}$ & $\varepsilon=10^{-8}$ \\
\hhline{|=||=|=|=|=|=|}
AMG-3G ($\boldsymbol{u}^{(0)}=\boldsymbol{0}$) & $\times$1 & $\times$1 & $\times$1 & $\times$1 & $\times$1\\ 
\hline
POD-2G ($\boldsymbol{u}^{(0)}=\boldsymbol{0}$) & 1.89 & 1.52 & 1.41 & 1.32 & 1.21\\
\hline
POD-2G ($\boldsymbol{u}^{(0)}=\boldsymbol{u}_{sur}$) & 12.77 & 6.10 & 3.06 & 2.19 & 1.76\\
\hline
\end{tabular}
\caption{Computational speedup of different preconditioners compared to the AMG-3G preconditioner}
\label{table:EX1_ComparisonPCG}
\end{table}
 
\noindent Again, the results obtained proved that the proposed methodology is superior than classic AMG preconditioners. In particular, for $\varepsilon = 10^{-4}$ and $\boldsymbol{u}^{(0)} = \boldsymbol{0}$, a reduction of computational cost of  $\times 1.89$ is observed between the proposed method and the 3-grid AMG. In addition, the initial solution delivered by the surrogate model, $\boldsymbol{u}^{(0)} = \boldsymbol{u}_{sur}$, is again a crucial factor of fast convergence, and can lead to a speedup of $\times 12.77$ for the same case.

\begin{figure}[H]
\centering
\begin{subfigure}[b]{0.49\textwidth}
         \centering
         \includegraphics[width=\textwidth]{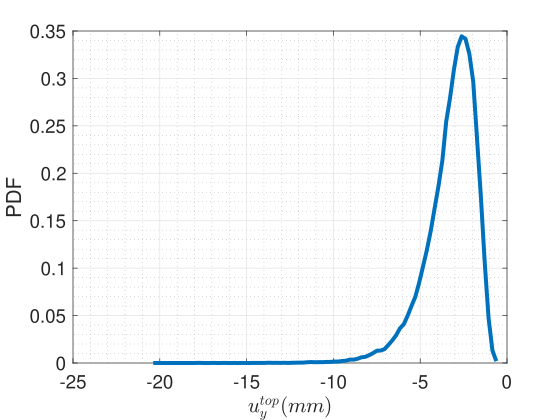}
         \caption{PDF of $u_{y}^{top}$}
         \label{fig:ex1_PDF}
     \end{subfigure}
     \hfill
     \begin{subfigure}[b]{0.49\textwidth}
         \centering
         \includegraphics[width=\textwidth]{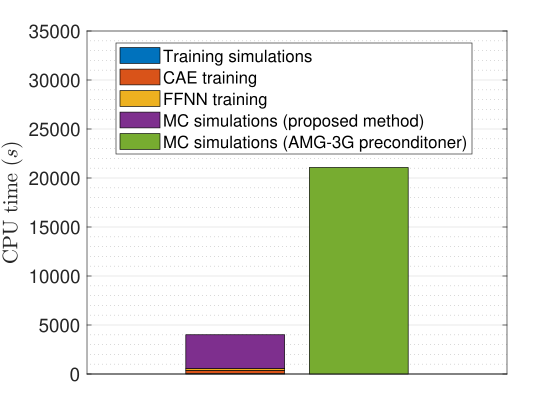}
          \caption{Comparison of computational cost}
          \label{fig:ex1_cost}
     \end{subfigure} 
         \caption{PDF of $u_{y}^{top}$ for $10^{5}$ MC simulations and comparison of computational cost.}
 \end{figure}

 Finally, in order to highlight the computational gain of the proposed framework in the context of the Monte Carlo method, $N_{MC} = 10^5$ simulations are performed to determine the probability density function (PDF) of the vertical displacement $u_{y}^{top}$ of the top node, where the load $P$ is applied. The calculated PDF is presented in figure \ref{fig:ex1_PDF}. Each simulation is solved with PCG and two different preconditioners, namely the proposed POD-2G method and a standard three grid Ruge-St{\"u}ben AMG preconditioner. The results are displayed in figure \ref{fig:ex1_cost} and verify that the proposed method is superior to classic AMG when dealing with parametrized systems. In particular, the conventional method needed 21109 $s$ to complete the $10^5$ simulations, while the proposed data-driven solver required 4013 $s$ for the same task including the offline cost (initial simulations and training of the surrogate model). This translates to a remarkable decrease in CPU time of $\times 5.26 $.

\subsection{Biot problem - deformable porous medium} 
\begin{figure}[H]
    \centering
    \includegraphics[width=0.5\textwidth]{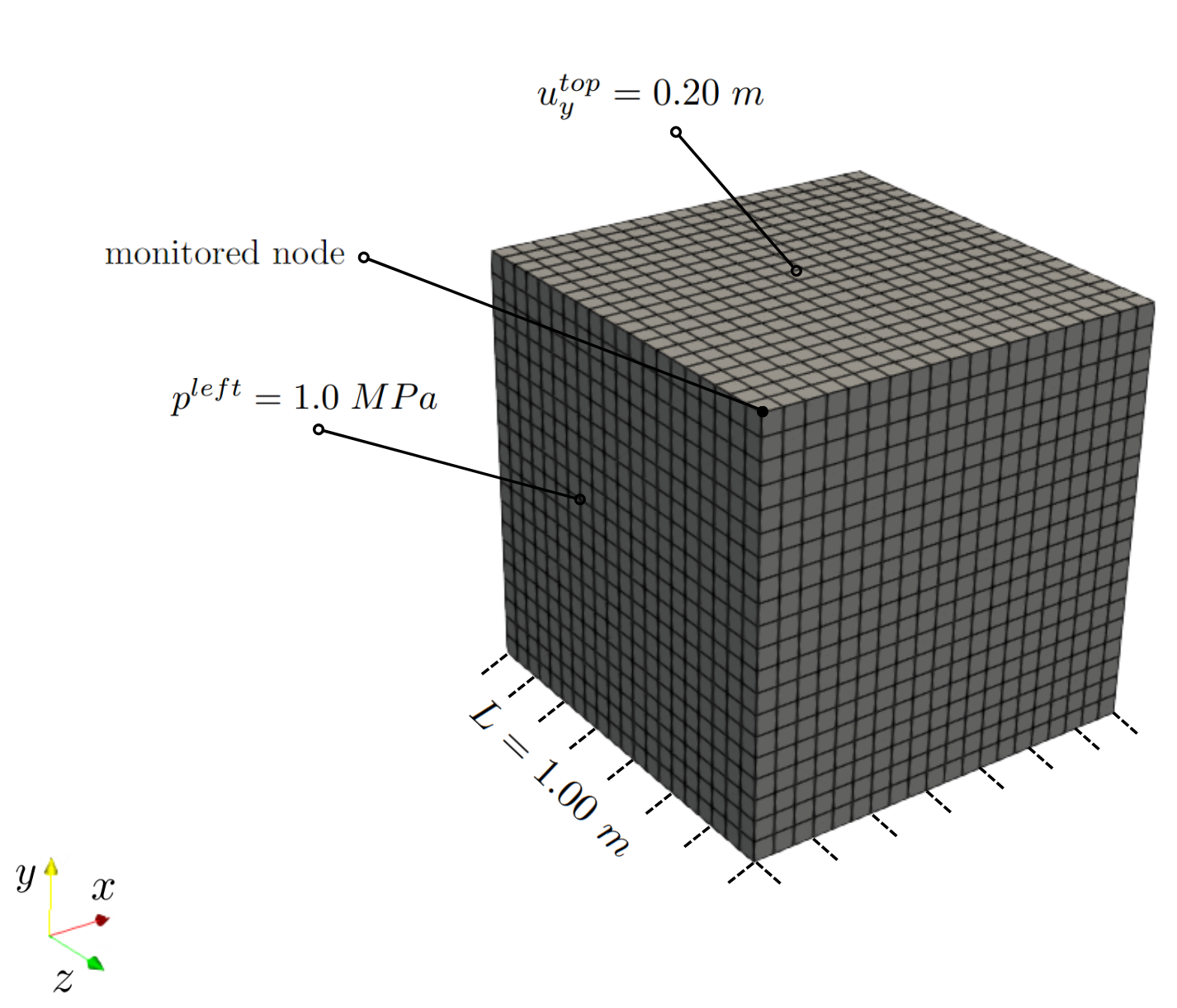}
    \caption{Geometry and boundary conditions of Biot problem}
    \label{fig:biot_bc}
\end{figure}

Biot's theory describes wave propagation in a porous saturated medium, i.e., a medium made of a solid matrix, fully soaked with a fluid. Biot does not take into account the microscopic level and assumes that continuum mechanics can be applied to measurable macroscopic quantities \cite{2001219}. Biot problem in weak form can be stated as:  Find $\boldsymbol{v}\in\mathcal{V}(\Omega;\mathbb{R}^3)$ and $p\in\mathcal{V}(\Omega;\mathbb{R})$ such that

\begin{equation}
\begin{aligned}
 \int_{\Omega} \boldsymbol{\sigma}(\boldsymbol{v}): \boldsymbol{\epsilon}(\boldsymbol{w})d\Omega
    - \int_{\Omega}  p\boldsymbol{A}: \boldsymbol{\epsilon}(\boldsymbol{w})d\Omega
    = 0
    &, \quad \forall \boldsymbol{w}\in\mathcal{V}_c(\Omega;\mathbb{R}^3)  \\
\int_{\Omega} q\boldsymbol{A}: \boldsymbol{\epsilon}(\boldsymbol{v})d\Omega + \int_{\Omega} \nabla q \cdot \boldsymbol{D}\ (\nabla p)^T d\Omega
    = 0
    &, \quad \forall q\in\mathcal{V}_c(\Omega;\mathbb{R}) \\
\boldsymbol{\sigma} = \lambda tr\left(\boldsymbol{\epsilon} \right)\mathbb{I}+2\mu\boldsymbol{\epsilon}&
\end{aligned} 
\label{eq:Biot}
\end{equation}  

\noindent with $\boldsymbol{A},\boldsymbol{D}$ being the Biot coefficient tensor and diffusion tensor, respectively. In this test case, the domain $\Omega$ is a cube and each side has a length of $L = 1.00 \ m$. Regarding the boundary conditions, a pressure distribution $p^{left}\coloneqq p|_{x=0} = 1.0 \ MPa$ is applied on the left face of the cube along with a displacement load $u_y^{top}\coloneqq u_y|_{z=1} = 0.20 \ m$ on the top face, while all displacements $u_x, u_y$ and $u_z$ are restrained in the bottom face ($z=0$). The problem definition is presented in figure \ref{fig:biot_bc}.

The finite element mesh contains 3-d hexa elements and the solution vector $\boldsymbol{u}\in{\mathbb{R}^d}$ consists of the nodal values of displacements and pressure, where in this case the total number of dofs is $d = 34839$. The Lame's constants $\mu$ and $\lambda$ are considered uncorrelated random variables following the lognormal distribution as described in table \ref{table:RandomParameters2}. The Poisson ratio $\nu$ is determined by: 
\begin{equation}
\nu = \frac{\lambda}{2(\lambda + \mu)} < 0.5
\end{equation} 

\noindent We further assumed that the Biot coefficient tensor $\boldsymbol{A}$ and $\boldsymbol{D}$ are constant, taking the values:

\begin{equation}
    \boldsymbol{A}=\begin{bmatrix} 0.13 & 0.13 & 0.13 \\ 0.09 & 0.09 & 0.09 \\ 0 & 0 & 0 \end{bmatrix}, \quad \boldsymbol{D}=\begin{bmatrix} 2.0 & 0.2 & 0 \\ 0.2 & 2.0 & 0 \\ 0 & 0 & 0.5 \end{bmatrix}
\end{equation}

\begin{table}[H]
 \centering
\begin{tabular}{|c|c|c|c|} 
\hline
 Parameter & Distribution & Mean & Standard deviation \\
\hline
$\mu(MPa)$ & Lognornal & $0.30$ & $0.09$\\ 
\hline
$\lambda(MPa)$ & Lognormal & $1.70$ & $0.51$\\ 
\hline
\end{tabular}
\caption{Random parameters of the Biot problem}
\label{table:RandomParameters2}
\end{table}

\begin{figure}[H]
\centering
\begin{subfigure}[b]{0.49\textwidth}
         \centering
         \includegraphics[width=\textwidth]{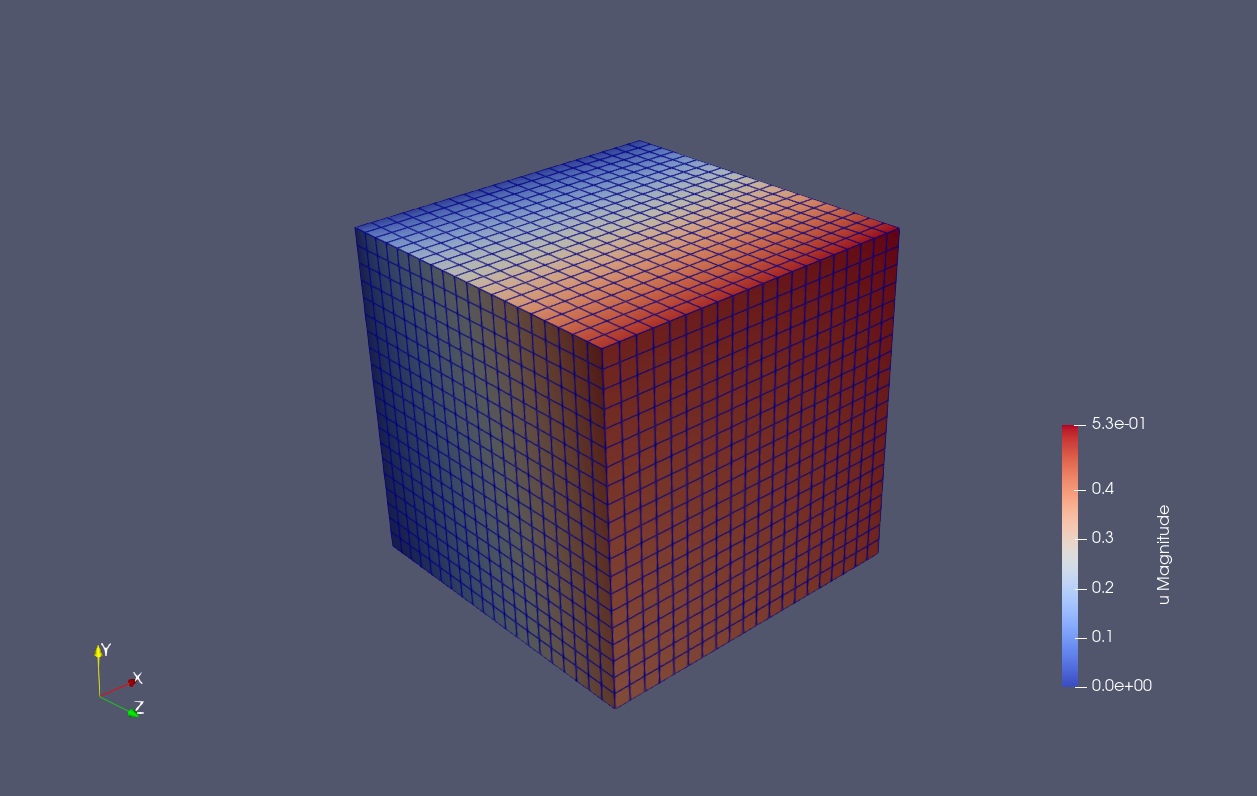}
         \caption{Displacement magnitude $\Vert\boldsymbol{u}\Vert$}
     \end{subfigure}
     \hfill
     \begin{subfigure}[b]{0.49\textwidth}
         \centering
         \includegraphics[width=\textwidth]{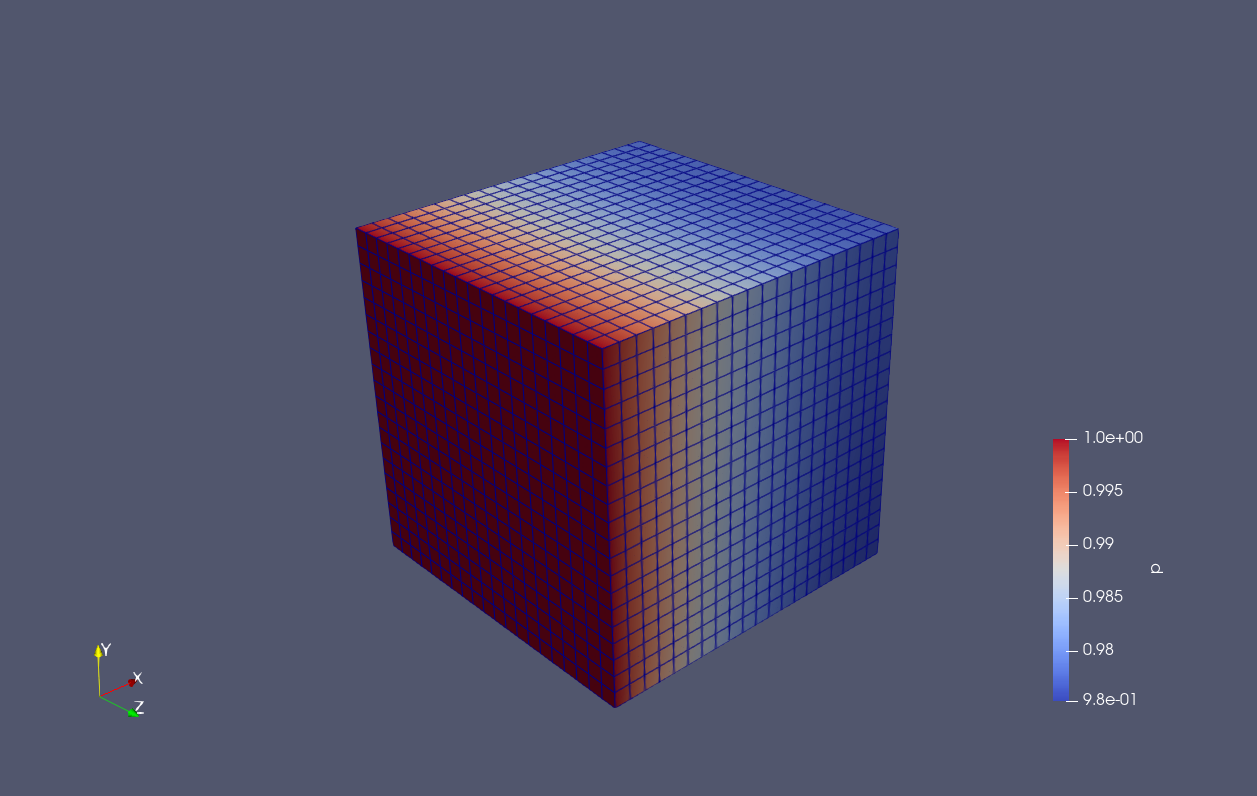}
          \caption{Pressure distribution $p$}
     \end{subfigure} 
         \caption{Displacement magnitude $\Vert \boldsymbol{u}\Vert$ and pressure distribution $p$ for $\lambda = 1.70 \ MPa$ and $\mu = 0.30 \ MPa$ }
    \label{fig:biot}
 \end{figure}

\noindent Figure \ref{fig:biot_bc} also displays a contour plot of the magnitude of $\boldsymbol{u}$ and the pressure distribution $p$ for $\mu = 0.30 \ MPa$ and $\lambda = 1.70\ MPa$.

\begin{figure}[H]
    \centering
    \includegraphics[width=1\textwidth]{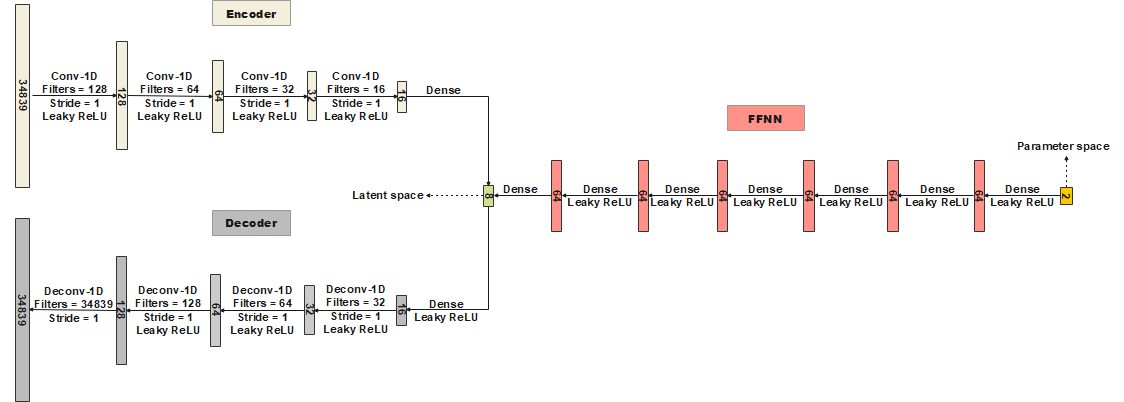}
    \caption{Surrogate model architecture}
    \label{fig:surrogate_2}
\end{figure}

The first step of the proposed methodology is to create an initial solution space. To this purpose, the Latin Hybercube sampling method was utilized to generate $N = 300$ parameter samples $\{[\mu_i, \lambda_i]\}_{i=1}^{N}$. The next steps are similar with those of the previous numerical example. The surrogate's architecture is presented in figure \ref{fig:surrogate_2}. The CAE is trained for 100 epochs with a batch size of 10 and a learning rate of $10^{-3}$, while the FFNN is trained for 5000 epochs with a batch size of 20 and a learning rate of $10^{-4}$. The average normalized $l_2$ norm error of the surrogate model in the test data set is 0.68\%. 

 As in the previous numerical example, a number of $N_{test} = 500$ parameter vectors $\{[\mu_j, \lambda_j]\}_{j=1}^{N_{test}}$ were generated according to their distribution and the corresponding problems were solved with the proposed POD-based solver and different Ruge-St{\"u}ben AMG solvers, with the number of grids ranging from 2 to 6. The size of the system of equations at the coarsest level for each of these solvers is presented in table \ref{table:ProblemSizeEx2}. For this example, 8 eigenvectors were retained in the POD expansion, as these were sufficient for capturing $99.99\%$ of the dataset's variance.
 
  \begin{table}[H]
 \centering
\begin{tabular}{|c|c|} 
\hline
 {} & System Size  \\
\hline
Initial Problem & $34839\times 34839$\\ 
\hline
AMG-2G & $8625 \times 8625$\\ 
\hline
AMG-3G & $1421 \times 1421$\\ 
\hline
AMG-4G & $229 \times 229$\\ 
\hline
AMG-5G & $47 \times 47$\\ 
\hline
AMG-6G & $9 \times 9$\\ 
\hline
POD-2G & $8 \times 8$\\ 
\hline
\end{tabular}
\caption{Size of the problem at the coarsest grid for the different solvers}
\label{table:ProblemSizeEx2}
\end{table}

\begin{figure}[H]
\centering
\begin{subfigure}[b]{0.49\textwidth}
         \centering
         \includegraphics[width=\textwidth]{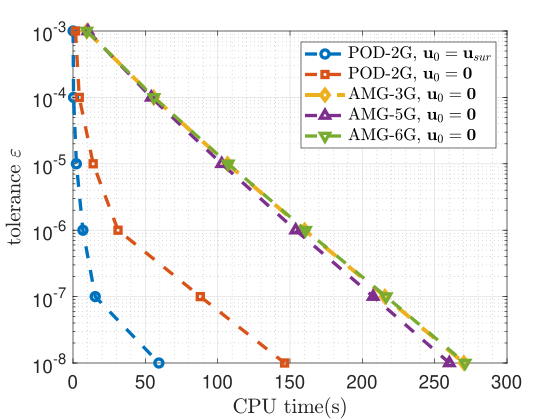}
         \caption{Comparison of mean CPU time}
     \end{subfigure}
     \hfill
     \begin{subfigure}[b]{0.49\textwidth}
         \centering
         \includegraphics[width=\textwidth]{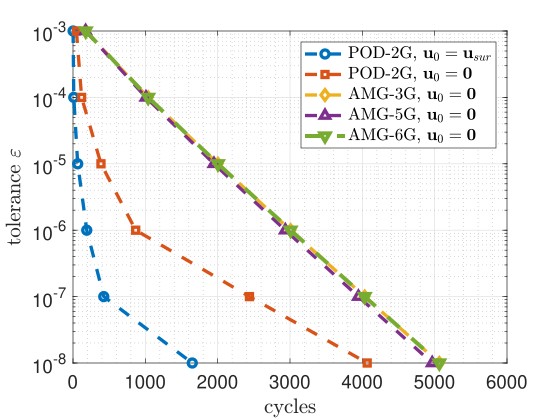}
          \caption{Comparison of mean number of cycles}
     \end{subfigure} 
         \caption{Comparison of mean CPU time and mean number of cycles over 500 analyses for different multigrid solvers}
    \label{fig:Results_AMG_2}
 \end{figure}

 The mean value of the CPU time and the number of cycles required for convergence to the desired number of tolerance are displayed in figure \ref{fig:Results_AMG_2} and table \ref{table:EX2_ComparisonAMG}. The results are very promising in terms of computational cost. For instance, for $\varepsilon = 10^{-5}$ and $\boldsymbol{u}^{(0)} = \boldsymbol{0}$, a reduction of computational cost of  $\times 7.32$ is achieved when comparing the proposed solver with the 3-grid AMG solver. Furthermore, obtaining an accurate initial solution $\boldsymbol{u}^{(0)}$ is again a very important component of the proposed framework. Specifically, by considering $\boldsymbol{u}^{(0)} = \boldsymbol{u}_{sur}$ instead of $\boldsymbol{u}^{(0)} = \boldsymbol{0}$ for $\varepsilon = 10^{-5}$, an additional decrease in CPU time of $\times 4.31$ can be achieved.

 \begin{table}[H]
 \centering
\begin{tabular}{|c||c|c|c|c|c|} 
\hline
  & $\varepsilon=10^{-4}$ & $\varepsilon=10^{-5}$ & $\varepsilon=10^{-6}$ & $\varepsilon=10^{-7}$ & $\varepsilon=10^{-8}$ \\
\hhline{|=||=|=|=|=|=|}
AMG-3G ($\boldsymbol{u}^{(0)}=\boldsymbol{0}$) & $\times$1 & $\times$1 & $\times$1 & $\times$1 & $\times$1\\ 
\hline
AMG-5G ($\boldsymbol{u}^{(0)}=\boldsymbol{0}$) & $\times$0.97 & $\times$0.96 & $\times$0.96 & $\times$0.96 & $\times$0.96\\ 
\hline
AMG-6G ($\boldsymbol{u}^{(0)}=\boldsymbol{0}$) & $\times$0.97 & $\times$0.96 & $\times$0.96 & $\times$0.96 & $\times$0.96\\ 
\hline
POD-2G ($\boldsymbol{u}^{(0)}=\boldsymbol{0}$) & $\times$12.31 & $\times$7.32 & $\times$4.89 & $\times$2.34 & $\times$1.77\\
\hline
POD-2G ($\boldsymbol{u}^{(0)}=\boldsymbol{u}_{sur}$) & $\times$76.89 & $\times$31.54 &$\times$17.90 & $\times$12.12 & $\times$4.35\\
\hline
\end{tabular}
\caption{Computational speedup of different solvers compared to AMG-3G}
\label{table:EX2_ComparisonAMG}
\end{table}
 
\begin{figure}[H]
\centering
\begin{subfigure}[b]{0.49\textwidth}
         \centering
         \includegraphics[width=\textwidth]{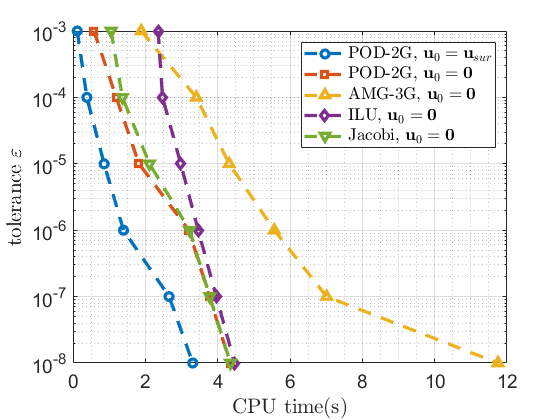}
         \caption{Comparison of mean CPU time}
     \end{subfigure}
     \hfill
     \begin{subfigure}[b]{0.49\textwidth}
         \centering
         \includegraphics[width=\textwidth]{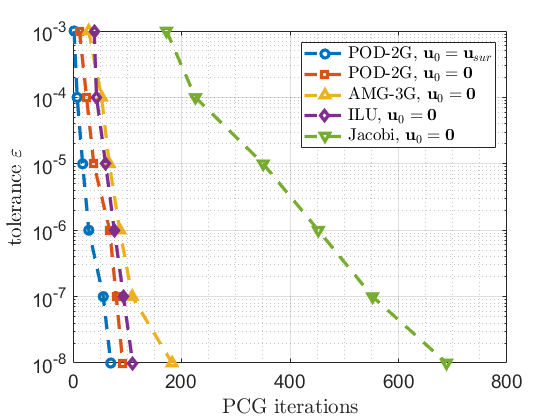}
          \caption{Comparison of mean PCG iterations}
     \end{subfigure} 
         \caption{Comparison of mean CPU time and mean number of PCG iterations over 500 analyses for different preconditioners}
    \label{fig:Results_PCG_2}
\end{figure}

 Furthermore, the convergence behaviour of the proposed method when used as a preconditioner in the context of the PCG method is presented in figure \ref{fig:Results_PCG_2}. Again, the results delivered by the proposed methodology showed its superior performance not only over AMG preconditioners but also over ILU and Jacobi preconditioners. In this case, for $\varepsilon = 10^{-5}$ and $\boldsymbol{u}^{(0)} = \boldsymbol{0}$, a reduction of computational cost of $\times 2.37$ is observed between the proposed method and the 3-grid AMG, of $\times 1.63$ with the ILU and of $\times 1.16$ with the Jacobi. Last but not least, the initial solution delivered by the surrogate model, $\boldsymbol{u}^{(0)} = \boldsymbol{u}_{sur}$, managed to further reduce the computational time by $\times 2.12$ when compared to POD-2G with $\boldsymbol{u}^{(0)} = \boldsymbol{0}$.

 \begin{table}[H]
 \centering
\begin{tabular}{|c||c|c|c|c|c|} 
\hline
  & $\varepsilon=10^{-4}$ & $\varepsilon=10^{-5}$ & $\varepsilon=10^{-6}$ & $\varepsilon=10^{-7}$ & $\varepsilon=10^{-8}$ \\
\hhline{|=||=|=|=|=|=|}
AMG-3G ($\boldsymbol{u}^{(0)}=\boldsymbol{0}$) & $\times$1 & $\times$1 & $\times$1 & $\times$1 & $\times$1\\
\hline
ILU ($\boldsymbol{u}^{(0)}=\boldsymbol{0}$) & $\times$1.38 & $\times$1.45 & $\times$1.61 & $\times$1.77 & $\times$2.63\\
\hline
Jacobi ($\boldsymbol{u}^{(0)}=\boldsymbol{0}$) & $\times$2.50 & $\times$2.04 & $\times$1.73 & $\times$1.85 & $\times$2.70\\
\hline$\times$
POD-2G ($\boldsymbol{u}^{(0)}=\boldsymbol{0}$) & $\times$2.86 & $\times$2.37 & $\times$1.74 & $\times$1.86 & $\times$2.71\\
\hline
POD-2G ($\boldsymbol{u}^{(0)}=\boldsymbol{u}_{sur}$) & $\times$8.88 & $\times$5.02 & $\times$3.98 & $\times$2.64 & $\times$3.55\\
\hline
\end{tabular}
\caption{Computational speedup of different preconditioners compared to the AMG-3G preconditioner}
\label{table:EX2_ComparisonPCG}
\end{table}

\begin{figure}[H]
\centering
\begin{subfigure}[b]{0.49\textwidth}
         \centering
         \includegraphics[width=\textwidth]{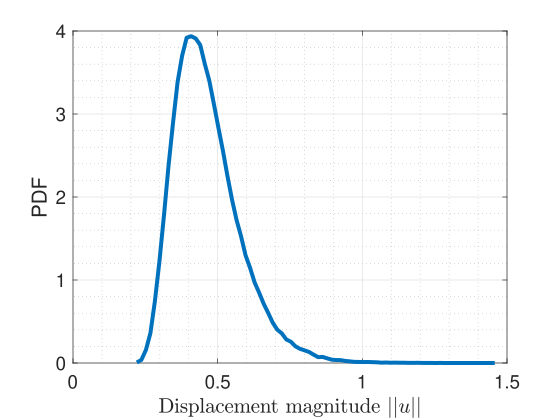}
         \caption{PDF of $\Vert \boldsymbol{u}\Vert$ at monitored dof}
         \label{fig:ex2_PDF}
     \end{subfigure}
     \hfill
     \begin{subfigure}[b]{0.49\textwidth}
         \centering
         \includegraphics[width=\textwidth]{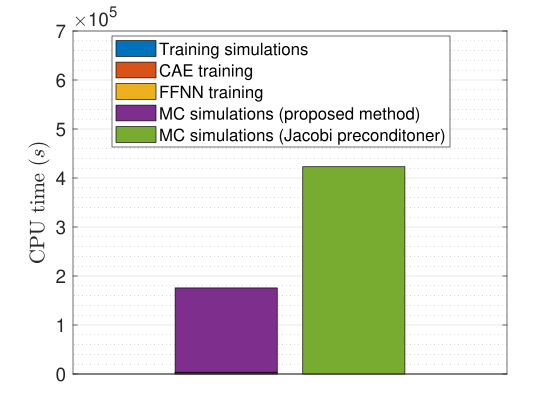}
          \caption{Comparison of computational cost}
          \label{fig:ex2_cost}
     \end{subfigure} 
         \caption{PDF of $\Vert \boldsymbol{u}\Vert$ at monitored dof for $2\times 10^{5}$ MC simulations and comparison of computational cost}
 \end{figure}

Finally, a Monte Carlo simulation is performed on this example as well, using $N_{MC} = 2\times 10^5$ simulations to determine the PDF of the displacement magnitude $\Vert \boldsymbol{u} \Vert$ of the monitored node (see figure \ref{fig:biot_bc}). The calculated PDF is presented in figure \ref{fig:ex2_PDF}. As in the previous example, each simulation is solved with PCG and two different preconditioners, namely the proposed POD-2G and a standard three grid Ruge-St{\"u}ben AMG preconditioner. Again, the results obtained by the proposed methods demonstrate a significant computational advantage over conventional preconditioners. In particular, the Jacobi preconditioner needed $4.23\times 10^{5}$ $s$ to complete $2 \times 10^{5}$ simulations, while the proposed data-driven solver required $1.75\times 10^{5}$ $s$ for the same task including the offline cost (initial simulations and training of the surrogate model). This translates to a decrease in CPU time of $\times 2.42$.

\section{Conclusions} \label{sec5}
The present work introduces a framework for accelerating the solution of parametrized problems that require multiple model evaluations. The proposed framework consists of two distinct yet complementary steps. The first step in the methodology is the construction of a `cheap-to-evaluate' metamodel using FFNNs and CAEs, trained over a reduced set of high-fidelity system solutions. Despite giving very accurate predictions at new parameter instances, these predictions are bound to exhibit some discrepancy with respect to the actual system solutions since they are not constrained by any physical laws. The second step in the methodology aims precisely at fixing this by proposing a data-driven iterative solver, inspired by the AMG method, that will refine the metamodel's predictions until a prescribed level of accuracy has been attained. In particular, using again the already available set of high-fidelity system solutions, POD is performed on this set to identify the subspace that captures most of the variation in the system responses. Next, a 2-level multigrid scheme is developed, termed POD-2G, using the projection operator from POD as the prolongation operator. This scheme was tested on numerical applications as a standalone solver, as well as a preconditioner to PCG, and in both cases, its superior performance with respect to conventional iterative solvers was demonstrated.

\bibliographystyle{elsarticle-num} 
\bibliography{main}

\begin{thebibliography}{10}
\expandafter\ifx\csname url\endcsname\relax
  \def\url#1{\texttt{#1}}\fi
\expandafter\ifx\csname urlprefix\endcsname\relax\def\urlprefix{URL }\fi
\expandafter\ifx\csname href\endcsname\relax
  \def\href#1#2{#2} \def\path#1{#1}\fi

\bibitem{Barrett1994}
R.~Barrett, M.~Berry, T.~F. Chan, J.~Demmel, J.~Donato, J.~Dongarra,
  V.~Eijkhout, R.~Pozo, C.~Romine, H.~V. der Vorst, Templates for the Solution
  of Linear Systems: Building Blocks for Iterative Methods, 2nd Edition, SIAM,
  Philadelphia, PA, 1994.

\bibitem{Benzi2000}
M.~Benzi, J.~K. Cullum, M.~Tuma, Robust approximate inverse preconditioning for
  the conjugate gradient method, SIAM Journal on Scientific Computing 22~(4)
  (2000) 1318--1332.
\newblock \href {https://doi.org/10.1137/S1064827599356900}
  {\path{doi:10.1137/S1064827599356900}}.

\bibitem{Lin2014}
F.-R. Lin, S.-W. Yang, X.-Q. Jin, Preconditioned iterative methods for
  fractional diffusion equation, Journal of Computational Physics 256 (2014)
  109 – 117.
\newblock \href {https://doi.org/10.1016/j.jcp.2013.07.040}
  {\path{doi:10.1016/j.jcp.2013.07.040}}.

\bibitem{Herzog2010}
R.~Herzog, E.~Sachs, Preconditioned conjugate gradient method for optimal
  control problems with control and state constraints, SIAM Journal on Matrix
  Analysis and Applications 31~(5) (2010) 2291--2317.
\newblock \href {https://doi.org/10.1137/090779127}
  {\path{doi:10.1137/090779127}}.

\bibitem{Saad1986}
Y.~Saad, M.~H. Schultz, Gmres: A generalized minimal residual algorithm for
  solving nonsymmetric linear systems, SIAM Journal on Scientific and
  Statistical Computing 7~(3) (1986) 856--869.
\newblock \href {https://doi.org/10.1137/0907058} {\path{doi:10.1137/0907058}}.

\bibitem{Shakib1989}
F.~Shakib, T.~J. Hughes, Z.~Johan, A multi-element group preconditioned gmres
  algorithm for nonsymmetric systems arising in finite element analysis,
  Computer Methods in Applied Mechanics and Engineering 75~(1-3) (1989) 415 –
  456.
\newblock \href {https://doi.org/10.1016/0045-7825(89)90040-6}
  {\path{doi:10.1016/0045-7825(89)90040-6}}.

\bibitem{Baglama1998}
J.~Baglama, D.~Calvetti, G.~Golub, L.~Reichel, Adaptively preconditioned gmres
  algorithms, SIAM Journal on Scientific Computing 20~(1) (1998) 243 – 269.
\newblock \href {https://doi.org/10.1137/S1064827596305258}
  {\path{doi:10.1137/S1064827596305258}}.

\bibitem{Concus1985}
P.~Concus, G.~H. Golub, G.~Meurant, Block preconditioning for the conjugate
  gradient method, SIAM Journal on Scientific and Statistical Computing 6~(1)
  (1985) 220--252.
\newblock \href {https://doi.org/10.1137/0906018} {\path{doi:10.1137/0906018}}.

\bibitem{Toselli2005}
A.~Toselli, O.~B. Widlund, Domain decomposition methods : algorithms and
  theory, 2005.

\bibitem{TALLEC1991}
P.~Tallec, Y.~Roeck, M.~Vidrascu,
  \href{https://www.sciencedirect.com/science/article/pii/037704279190150I}{Domain
  decomposition methods for large linearly elliptic three-dimensional
  problems}, Journal of Computational and Applied Mathematics 34~(1) (1991)
  93--117.
\newblock \href {https://doi.org/https://doi.org/10.1016/0377-0427(91)90150-I}
  {\path{doi:https://doi.org/10.1016/0377-0427(91)90150-I}}.
\newline\urlprefix\url{https://www.sciencedirect.com/science/article/pii/037704279190150I}

\bibitem{Farhat1991}
C.~Farhat, F.-X. Roux, A method of finite element tearing and interconnecting
  and its parallel solution algorithm, International Journal for Numerical
  Methods in Engineering 32~(6) (1991) 1205--1227.
\newblock \href {https://doi.org/https://doi.org/10.1002/nme.1620320604}
  {\path{doi:https://doi.org/10.1002/nme.1620320604}}.

\bibitem{FARHAT1994}
C.~Farhat, J.~Mandel, F.~X. Roux, Optimal convergence properties of the feti
  domain decomposition method, Computer Methods in Applied Mechanics and
  Engineering 115~(3) (1994) 365--385.
\newblock \href {https://doi.org/https://doi.org/10.1016/0045-7825(94)90068-X}
  {\path{doi:https://doi.org/10.1016/0045-7825(94)90068-X}}.

\bibitem{Fragakis2003}
Y.~Fragakis, M.~Papadrakakis, The mosaic of high performance domain
  decomposition methods for structural mechanics: Formulation, interrelation
  and numerical efficiency of primal and dual methods, Computer Methods in
  Applied Mechanics and Engineering 192 (2003) 3799--3830.

\bibitem{Cai1999}
X.-C. Cai, M.~Sarkis, A restricted additive schwarz preconditioner for general
  sparse linear systems, SIAM Journal on Scientific Computing 21~(2) (1999)
  792--797.
\newblock \href {https://doi.org/10.1137/S106482759732678X}
  {\path{doi:10.1137/S106482759732678X}}.

\bibitem{DAAS2021}
H.~Al~Daas, L.~Grigori, P.~Jolivet, P.-H. Tournier, A multilevel schwarz
  preconditioner based on a hierarchy of robust coarse spaces, SIAM Journal on
  Scientific Computing 43~(3) (2021) A1907--A1928.
\newblock \href {https://doi.org/10.1137/19M1266964}
  {\path{doi:10.1137/19M1266964}}.

\bibitem{Trottenberg2000}
U.Trottenberg, C.~Oostelee, A.~Schuller, Multigrid 1st Edition, Academic Press,
  2000.

\bibitem{Iwamura2003}
C.~Iwamura, F.~S. Costa, I.~Sbarski, A.~Easton, N.~Li, An efficient algebraic
  multigrid preconditioned conjugate gradient solver, Computer Methods in
  Applied Mechanics and Engineering 192~(20-21) (2003) 2299 – 2318.
\newblock \href {https://doi.org/10.1016/S0045-7825(02)00378-X}
  {\path{doi:10.1016/S0045-7825(02)00378-X}}.

\bibitem{Heys2005}
J.~Heys, T.~Manteuffel, S.~F. McCormick, L.~Olson, Algebraic multigrid for
  higher-order finite elements, Journal of Computational Physics 204~(2) (2005)
  520 – 532.
\newblock \href {https://doi.org/10.1016/j.jcp.2004.10.021}
  {\path{doi:10.1016/j.jcp.2004.10.021}}.

\bibitem{Langer2003}
U.~Langer, D.~Pusch, S.~Reitzinger, Efficient preconditioners for boundary
  element matrices based on grey-box algebraic multigrid methods, International
  Journal for Numerical Methods in Engineering 58~(13) (2003) 1937 – 1953.
\newblock \href {https://doi.org/10.1002/nme.839} {\path{doi:10.1002/nme.839}}.

\bibitem{RAMAGE1999}
A.~Ramage,
  \href{https://www.sciencedirect.com/science/article/pii/S0377042799002344}{A
  multigrid preconditioner for stabilised discretisations of
  advection–diffusion problems}, Journal of Computational and Applied
  Mathematics 110~(1) (1999) 187--203.
\newblock \href {https://doi.org/https://doi.org/10.1016/S0377-0427(99)00234-4}
  {\path{doi:https://doi.org/10.1016/S0377-0427(99)00234-4}}.
\newline\urlprefix\url{https://www.sciencedirect.com/science/article/pii/S0377042799002344}

\bibitem{Wienands2000}
R.~Wienands, C.~W. Oosterlee, T.~Washio, Fourier analysis of gmres(m)
  preconditioned by multigrid, SIAM Journal on Scientific Computing 22~(2)
  (2000) 582--603.
\newblock \href {https://doi.org/10.1137/S1064827599353014}
  {\path{doi:10.1137/S1064827599353014}}.

\bibitem{Vakili2009}
S.~Vakili, M.~Darbandi, Recommendations on enhancing the efficiency of
  algebraic multigrid preconditioned gmres in solving coupled fluid flow
  equations, Numerical Heat Transfer, Part B: Fundamentals 55~(3) (2009) 232
  – 256.
\newblock \href {https://doi.org/10.1080/10407790802628879}
  {\path{doi:10.1080/10407790802628879}}.

\bibitem{Zahm2016}
O.~Zahm, A.~Nouy, Interpolation of inverse operators for preconditioning
  parameter-dependent equations, SIAM Journal on Scientific Computing 38~(2)
  (2016) A1044--A1074.
\newblock \href {https://doi.org/10.1137/15M1019210}
  {\path{doi:10.1137/15M1019210}}.

\bibitem{Bergamaschi2020}
L.~Bergamaschi, A survey of low-rank updates of preconditioners for sequences
  of symmetric linear systems, Algorithms 13~(4) (2020).
\newblock \href {https://doi.org/10.3390/A13040100}
  {\path{doi:10.3390/A13040100}}.

\bibitem{Carr2021}
A.~Carr, E.~de~Sturler, S.~Gugercin, Preconditioning parametrized linear
  systems, SIAM Journal on Scientific Computing 43~(3) (2021) A2242--A2267.
\newblock \href {https://doi.org/10.1137/20M1331123}
  {\path{doi:10.1137/20M1331123}}.

\bibitem{STAVROULAKIS2014}
G.~Stavroulakis, D.~G. Giovanis, M.~Papadrakakis, V.~Papadopoulos,
  \href{https://www.sciencedirect.com/science/article/pii/S0045782514000954}{A
  new perspective on the solution of uncertainty quantification and reliability
  analysis of large-scale problems}, Computer Methods in Applied Mechanics and
  Engineering 276 (2014) 627--658.
\newblock \href {https://doi.org/https://doi.org/10.1016/j.cma.2014.03.009}
  {\path{doi:https://doi.org/10.1016/j.cma.2014.03.009}}.
\newline\urlprefix\url{https://www.sciencedirect.com/science/article/pii/S0045782514000954}

\bibitem{Saad1997}
Y.~Saad, Analysis of augmented krylov subspace methods, SIAM Journal on Matrix
  Analysis and Applications 18~(2) (1997) 435--449.
\newblock \href {https://doi.org/10.1137/S0895479895294289}
  {\path{doi:10.1137/S0895479895294289}}.

\bibitem{Sturler1999}
E.~de~Sturler, Truncation strategies for optimal krylov subspace methods, SIAM
  Journal on Numerical Analysis 36~(3) (1999) 864--889.
\newblock \href {https://doi.org/10.1137/S0036142997315950}
  {\path{doi:10.1137/S0036142997315950}}.

\bibitem{Chapman1997}
A.~Chapman, Y.~Saad, Deflated and augmented krylov subspace techniques,
  Numerical Linear Algebra with Applications 4~(1) (1997) 43--66.
\newblock \href
  {https://doi.org/https://doi.org/10.1002/(SICI)1099-1506(199701/02)4:1<43::AID-NLA99>3.0.CO;2-Z}
  {\path{doi:https://doi.org/10.1002/(SICI)1099-1506(199701/02)4:1<43::AID-NLA99>3.0.CO;2-Z}}.

\bibitem{Saad2000}
Y.~Saad, M.~Yeung, J.~Erhel, F.~Guyomarc'h, A deflated version of the conjugate
  gradient algorithm, SIAM Journal on Scientific Computing 21~(5) (2000)
  1909--1926.
\newblock \href {https://doi.org/10.1137/S1064829598339761}
  {\path{doi:10.1137/S1064829598339761}}.

\bibitem{Daas2021a}
D.~H.A., G.~L., H.~P., R.~P., Recycling krylov subspaces and truncating
  deflation subspaces for solving sequence of linear systems, ACM Transactions
  on Mathematical Software 47~(2) (2021).
\newblock \href {https://doi.org/10.1145/3439746} {\path{doi:10.1145/3439746}}.

\bibitem{Papadrakakis1996}
M.~Papadrakakis, V.~Papadopoulos, N.~D. Lagaros, Structural reliability analyis
  of elastic-plastic structures using neural networks and monte carlo
  simulation, Computer Methods in Applied Mechanics and Engineering 136~(1-2)
  (1996) 145 – 163.
\newblock \href {https://doi.org/10.1016/0045-7825(96)01011-0}
  {\path{doi:10.1016/0045-7825(96)01011-0}}.

\bibitem{PAPADRAKAKIS2002}
M.~Papadrakakis, N.~D. Lagaros,
  \href{https://www.sciencedirect.com/science/article/pii/S0045782502002876}{Reliability-based
  structural optimization using neural networks and monte carlo simulation},
  Computer Methods in Applied Mechanics and Engineering 191~(32) (2002)
  3491--3507.
\newblock \href {https://doi.org/https://doi.org/10.1016/S0045-7825(02)00287-6}
  {\path{doi:https://doi.org/10.1016/S0045-7825(02)00287-6}}.
\newline\urlprefix\url{https://www.sciencedirect.com/science/article/pii/S0045782502002876}

\bibitem{SEYHAN2005}
A.~T. Seyhan, G.~Tayfur, M.~Karakurt, M.~Tanogˇlu,
  \href{https://www.sciencedirect.com/science/article/pii/S0927025604003076}{Artificial
  neural network (ann) prediction of compressive strength of vartm processed
  polymer composites}, Computational Materials Science 34~(1) (2005) 99--105.
\newblock \href
  {https://doi.org/https://doi.org/10.1016/j.commatsci.2004.11.001}
  {\path{doi:https://doi.org/10.1016/j.commatsci.2004.11.001}}.
\newline\urlprefix\url{https://www.sciencedirect.com/science/article/pii/S0927025604003076}

\bibitem{HOSNIELHEWY2006}
A.~{Hosni Elhewy}, E.~Mesbahi, Y.~Pu, Reliability analysis of structures using
  neural network method, Probabilistic Engineering Mechanics 21~(1) (2006)
  44--53.
\newblock \href
  {https://doi.org/https://doi.org/10.1016/j.probengmech.2005.07.002}
  {\path{doi:https://doi.org/10.1016/j.probengmech.2005.07.002}}.

\bibitem{Chojaczyk2015}
A.~Chojaczyk, A.~Teixeira, L.~Neves, J.~Cardoso, C.~Guedes~Soares, Review and
  application of artificial neural networks models in reliability analysis of
  steel structures, Structural Safety 52~(PA) (2015) 78 – 89.
\newblock \href {https://doi.org/10.1016/j.strusafe.2014.09.002}
  {\path{doi:10.1016/j.strusafe.2014.09.002}}.

\bibitem{NIKOLOPOULOS2021}
S.~Nikolopoulos, I.~Kalogeris, V.~Papadopoulos,
  \href{https://www.sciencedirect.com/science/article/pii/S0952197621004541}{Non-intrusive
  surrogate modeling for parametrized time-dependent partial differential
  equations using convolutional autoencoders}, Engineering Applications of
  Artificial Intelligence 109 (2022) 104652.
\newblock \href
  {https://doi.org/https://doi.org/10.1016/j.engappai.2021.104652}
  {\path{doi:https://doi.org/10.1016/j.engappai.2021.104652}}.
\newline\urlprefix\url{https://www.sciencedirect.com/science/article/pii/S0952197621004541}

\bibitem{NIKOLOPOULOS2022}
S.~Nikolopoulos, I.~Kalogeris, V.~Papadopoulos, Machine learning accelerated
  transient analysis of stochastic nonlinear structures, Engineering Structures
  257 (2022) 114020.
\newblock \href
  {https://doi.org/https://doi.org/10.1016/j.engstruct.2022.114020}
  {\path{doi:https://doi.org/10.1016/j.engstruct.2022.114020}}.

\bibitem{XU2020}
J.~Xu, K.~Duraisamy,
  \href{https://www.sciencedirect.com/science/article/pii/S0045782520305648}{Multi-level
  convolutional autoencoder networks for parametric prediction of
  spatio-temporal dynamics}, Computer Methods in Applied Mechanics and
  Engineering 372 (2020) 113379.
\newblock \href {https://doi.org/https://doi.org/10.1016/j.cma.2020.113379}
  {\path{doi:https://doi.org/10.1016/j.cma.2020.113379}}.
\newline\urlprefix\url{https://www.sciencedirect.com/science/article/pii/S0045782520305648}

\bibitem{Yu2019}
Y.~Y., Y.~H., L.~Y., Aircraft dynamics simulation using a novel physics-based
  learning method, Aerospace Science and Technology 87 (2019) 254 – 264.
\newblock \href {https://doi.org/10.1016/j.ast.2019.02.021}
  {\path{doi:10.1016/j.ast.2019.02.021}}.

\bibitem{Zhou2019}
Z.~J.M., D.~L., G.~W., Y.~J., Impact load identification of nonlinear
  structures using deep recurrent neural network, Mechanical Systems and Signal
  Processing 133 (2019).
\newblock \href {https://doi.org/10.1016/j.ymssp.2019.106292}
  {\path{doi:10.1016/j.ymssp.2019.106292}}.

\bibitem{Carlberg2011}
K.~Carlberg, C.~Farhat, A low-cost, goal-oriented ‘compact proper orthogonal
  decomposition’ basis for model reduction of static systems, International
  Journal for Numerical Methods in Engineering 86~(3) (2011) 381--402.
\newblock \href {https://doi.org/https://doi.org/10.1002/nme.3074}
  {\path{doi:https://doi.org/10.1002/nme.3074}}.

\bibitem{Zahr2017}
M.~J. Zahr, P.~Avery, C.~Farhat, A multilevel projection-based model order
  reduction framework for nonlinear dynamic multiscale problems in structural
  and solid mechanics, International Journal for Numerical Methods in
  Engineering 112~(8) (2017) 855--881.
\newblock \href {https://doi.org/https://doi.org/10.1002/nme.5535}
  {\path{doi:https://doi.org/10.1002/nme.5535}}.

\bibitem{Agathos2020}
K.~Agathos, S.~P.~A. Bordas, E.~Chatzi, Parametrized reduced order modeling for
  cracked solids, International Journal for Numerical Methods in Engineering
  121~(20) (2020) 4537 – 4565.
\newblock \href {https://doi.org/10.1002/nme.6447}
  {\path{doi:10.1002/nme.6447}}.

\bibitem{Chinesta2010}
F.~Chinesta, A.~Ammar, E.~Cueto, Proper generalized decomposition of multiscale
  models, International Journal for Numerical Methods in Engineering 83~(8-9)
  (2010) 1114--1132.
\newblock \href {https://doi.org/https://doi.org/10.1002/nme.2794}
  {\path{doi:https://doi.org/10.1002/nme.2794}}.

\bibitem{Ladeveze2010}
P.~Ladevèze, J.-C. Passieux, D.~Néron, The latin multiscale computational
  method and the proper generalized decomposition, Computer Methods in Applied
  Mechanics and Engineering 199~(21-22) (2010) 1287 – 1296.
\newblock \href {https://doi.org/10.1016/j.cma.2009.06.023}
  {\path{doi:10.1016/j.cma.2009.06.023}}.

\bibitem{Ladeveze2011}
P.~Ladevèze, L.~Chamoin, On the verification of model reduction methods based
  on the proper generalized decomposition, Computer Methods in Applied
  Mechanics and Engineering 200~(23-24) (2011) 2032 – 2047.
\newblock \href {https://doi.org/10.1016/j.cma.2011.02.019}
  {\path{doi:10.1016/j.cma.2011.02.019}}.

\bibitem{DALSANTO2020}
N.~{Dal Santo}, S.~Deparis, L.~Pegolotti,
  \href{https://www.sciencedirect.com/science/article/pii/S0021999120303247}{Data
  driven approximation of parametrized pdes by reduced basis and neural
  networks}, Journal of Computational Physics 416 (2020) 109550.
\newblock \href {https://doi.org/https://doi.org/10.1016/j.jcp.2020.109550}
  {\path{doi:https://doi.org/10.1016/j.jcp.2020.109550}}.
\newline\urlprefix\url{https://www.sciencedirect.com/science/article/pii/S0021999120303247}

\bibitem{SALVADOR2021}
M.~Salvador, L.~Dedè, A.~Manzoni,
  \href{https://www.sciencedirect.com/science/article/pii/S0898122121003928}{Non
  intrusive reduced order modeling of parametrized pdes by kernel pod and
  neural networks}, Computers {\&} Mathematics with Applications 104 (2021)
  1--13.
\newblock \href {https://doi.org/https://doi.org/10.1016/j.camwa.2021.11.001}
  {\path{doi:https://doi.org/10.1016/j.camwa.2021.11.001}}.
\newline\urlprefix\url{https://www.sciencedirect.com/science/article/pii/S0898122121003928}

\bibitem{KALOGERIS2021}
I.~Kalogeris, V.~Papadopoulos,
  \href{https://www.sciencedirect.com/science/article/pii/S0045782520307532}{Diffusion
  maps-aided neural networks for the solution of parametrized pdes}, Computer
  Methods in Applied Mechanics and Engineering 376 (2021) 113568.
\newblock \href {https://doi.org/https://doi.org/10.1016/j.cma.2020.113568}
  {\path{doi:https://doi.org/10.1016/j.cma.2020.113568}}.
\newline\urlprefix\url{https://www.sciencedirect.com/science/article/pii/S0045782520307532}

\bibitem{dosSantos2022}
K.~R. dos Santos, D.~G. Giovanis, M.~D. Shields, Grassmannian diffusion
  maps--based dimension reduction and classification for high-dimensional data,
  SIAM Journal on Scientific Computing 44~(2) (2022) B250--B274.
\newblock \href {https://doi.org/10.1137/20M137001X}
  {\path{doi:10.1137/20M137001X}}.

\bibitem{Kadeethum2022}
T.~Kadeethum, F.~Ballarin, Y.~Choi, D.~O'Malley, H.~Yoon, N.~Bouklas,
  Non-intrusive reduced order modeling of natural convection in porous media
  using convolutional autoencoders: Comparison with linear subspace techniques,
  Advances in Water Resources 160 (2022).
\newblock \href {https://doi.org/10.1016/j.advwatres.2021.104098}
  {\path{doi:10.1016/j.advwatres.2021.104098}}.

\bibitem{VLACHAS2021}
K.~Vlachas, K.~Tatsis, K.~Agathos, A.~R. Brink, E.~Chatzi,
  \href{https://www.sciencedirect.com/science/article/pii/S0022460X21001279}{A
  local basis approximation approach for nonlinear parametric model order
  reduction}, Journal of Sound and Vibration 502 (2021) 116055.
\newblock \href {https://doi.org/https://doi.org/10.1016/j.jsv.2021.116055}
  {\path{doi:https://doi.org/10.1016/j.jsv.2021.116055}}.
\newline\urlprefix\url{https://www.sciencedirect.com/science/article/pii/S0022460X21001279}

\bibitem{Heaney}
C.~E. Heaney, Z.~Wolffs, J.~A. Tómasson, L.~Kahouadji, P.~Salinas, A.~Nicolle,
  I.~M. Navon, O.~K. Matar, N.~Srinil, C.~C. Pain,
  \href{https://doi.org/10.1063/5.0088070}{An ai-based non-intrusive
  reduced-order model for extended domains applied to multiphase flow in
  pipes}, Physics of Fluids 34~(5) (2022) 055111.
\newblock \href {http://arxiv.org/abs/https://doi.org/10.1063/5.0088070}
  {\path{arXiv:https://doi.org/10.1063/5.0088070}}, \href
  {https://doi.org/10.1063/5.0088070} {\path{doi:10.1063/5.0088070}}.
\newline\urlprefix\url{https://doi.org/10.1063/5.0088070}

\bibitem{LAZZARA2022107629}
M.~Lazzara, M.~Chevalier, M.~Colombo, J.~{Garay Garcia}, C.~Lapeyre, O.~Teste,
  \href{https://www.sciencedirect.com/science/article/pii/S1270963822003030}{Surrogate
  modelling for an aircraft dynamic landing loads simulation using an lstm
  autoencoder-based dimensionality reduction approach}, Aerospace Science and
  Technology (2022) 107629\href
  {https://doi.org/https://doi.org/10.1016/j.ast.2022.107629}
  {\path{doi:https://doi.org/10.1016/j.ast.2022.107629}}.
\newline\urlprefix\url{https://www.sciencedirect.com/science/article/pii/S1270963822003030}

\bibitem{Carlberg2016}
K.~Carlberg, V.~Forstall, R.~Tuminaro, Krylov-subspace recycling via the
  pod-augmented conjugate-gradient method, SIAM Journal on Matrix Analysis and
  Applications 37~(3) (2016) 1304--1336.
\newblock \href {https://doi.org/10.1137/16M1057693}
  {\path{doi:10.1137/16M1057693}}.

\bibitem{Heinlein2019}
A.~Heinlein, A.~Klawonn, M.~Lanser, J.~Weber, Machine learning in adaptive
  domain decomposition methods---predicting the geometric location of
  constraints, SIAM Journal on Scientific Computing 41~(6) (2019) A3887--A3912.
\newblock \href {https://doi.org/10.1137/18M1205364}
  {\path{doi:10.1137/18M1205364}}.

\bibitem{CHEN2022}
Y.~Chen, B.~Dong, J.~Xu, Meta-mgnet: Meta multigrid networks for solving
  parameterized partial differential equations, Journal of Computational
  Physics 455 (2022) 110996.
\newblock \href {https://doi.org/https://doi.org/10.1016/j.jcp.2022.110996}
  {\path{doi:https://doi.org/10.1016/j.jcp.2022.110996}}.

\bibitem{luz2020}
I.~Luz, M.~Galun, H.~Maron, R.~Basri, I.~Yavneh, Learning algebraic multigrid
  using graph neural networks, in: International Conference on Machine
  Learning, PMLR, 2020, pp. 6489--6499.

\bibitem{Hestenes1952}
M.~R. Hestenes, E.~Stiefel, Methods of conjugate gradients for solving linear
  systems, Journal of research of the National Bureau of Standards 49 (1952)
  409--435.

\bibitem{Ruge1987}
J.~W. Ruge, K.~Stüben,
  \href{https://epubs.siam.org/doi/abs/10.1137/1.9781611971057.ch4}{Algebraic
  Multigrid}, Frontiers in Applied Mathematics, SIAM, 1987, Ch. 4. Algebraic
  Multigrid, pp. 73--130.
\newblock \href {https://doi.org/10.1137/1.9781611971057.ch4}
  {\path{doi:10.1137/1.9781611971057.ch4}}.
\newline\urlprefix\url{https://epubs.siam.org/doi/abs/10.1137/1.9781611971057.ch4}

\bibitem{Stuben2001}
K.~Stüben,
  \href{https://www.sciencedirect.com/science/article/pii/S0377042700005161}{A
  review of algebraic multigrid}, Journal of Computational and Applied
  Mathematics 128~(1) (2001) 281--309, numerical Analysis 2000. Vol. VII:
  Partial Differential Equations.
\newblock \href {https://doi.org/https://doi.org/10.1016/S0377-0427(00)00516-1}
  {\path{doi:https://doi.org/10.1016/S0377-0427(00)00516-1}}.
\newline\urlprefix\url{https://www.sciencedirect.com/science/article/pii/S0377042700005161}

\bibitem{Brezina2001}
M.~Brezina, A.~J. Cleary, R.~D. Falgout, V.~E. Henson, J.~E. Jones, T.~A.
  Manteuffel, S.~F. McCormick, J.~W. Ruge, Algebraic multigrid based on element
  interpolation (amge), SIAM Journal on Scientific Computing 22~(5) (2001)
  1570--1592.
\newblock \href {https://doi.org/10.1137/S1064827598344303}
  {\path{doi:10.1137/S1064827598344303}}.

\bibitem{Treister2010}
E.~Treister, I.~Yavneh, Square and stretch multigrid for stochastic matrix
  eigenproblems, Numerical Linear Algebra with Applications 17~(2-3) (2010)
  229--251.
\newblock \href {https://doi.org/https://doi.org/10.1002/nla.708}
  {\path{doi:https://doi.org/10.1002/nla.708}}.

\bibitem{Napov2016}
A.~Napov, Y.~Notay,
  \href{https://www.scopus.com/inward/record.uri?eid=2-s2.0-85010702842&partnerID=40&md5=41629148826b0c387598c23a886a6f81}{An
  efficient multigrid method for graph laplacian systems}, Electronic
  Transactions on Numerical Analysis 45 (2016) 201 – 218, cited by: 10.
\newline\urlprefix\url{https://www.scopus.com/inward/record.uri?eid=2-s2.0-85010702842&partnerID=40&md5=41629148826b0c387598c23a886a6f81}

\bibitem{Facca2021}
E.~Facca, M.~Benzi, Fast iterative solution of the optimal transport problem on
  graphs, SIAM Journal on Scientific Computing 43~(3) (2021) A2295 – A2319.
\newblock \href {https://doi.org/10.1137/20M137015X}
  {\path{doi:10.1137/20M137015X}}.

\bibitem{Vanek1996}
V.~P., M.~J., B.~M., Algrebraic multigrid by smoothed aggregation for second
  and fourth order elliptic problems, Computing 56~(3) (1996) 179--196.
\newblock \href {https://doi.org/10.1007/BF02238511}
  {\path{doi:10.1007/BF02238511}}.

\bibitem{Mo2019}
S.~Mo, Y.~Zhu, N.~Zabaras, X.~Shi, J.~Wu,
  \href{https://agupubs.onlinelibrary.wiley.com/doi/abs/10.1029/2018WR023528}{Deep
  convolutional encoder-decoder networks for uncertainty quantification of
  dynamic multiphase flow in heterogeneous media}, Water Resources Research
  55~(1) (2019) 703--728.
\newblock \href {https://doi.org/https://doi.org/10.1029/2018WR023528}
  {\path{doi:https://doi.org/10.1029/2018WR023528}}.
\newline\urlprefix\url{https://agupubs.onlinelibrary.wiley.com/doi/abs/10.1029/2018WR023528}

\bibitem{olsson2002latin}
A.~M. Olsson, G.~E. Sandberg, Latin hypercube sampling for stochastic finite
  element analysis, Journal of Engineering Mechanics 128~(1) (2002) 121--125.

\bibitem{Rathinam2003}
M.~Rathinam, L.~R. Petzold, A new look at proper orthogonal decomposition, SIAM
  Journal on Numerical Analysis 41~(5) (2003) 1893--1925.
\newblock \href {https://doi.org/10.1137/S0036142901389049}
  {\path{doi:10.1137/S0036142901389049}}.

\bibitem{LIEU2006}
T.~Lieu, C.~Farhat, M.~Lesoinne,
  \href{https://www.sciencedirect.com/science/article/pii/S0045782505005153}{Reduced-order
  fluid/structure modeling of a complete aircraft configuration}, Computer
  Methods in Applied Mechanics and Engineering 195~(41) (2006) 5730--5742, john
  H. Argyris Memorial Issue. Part II.
\newblock \href {https://doi.org/https://doi.org/10.1016/j.cma.2005.08.026}
  {\path{doi:https://doi.org/10.1016/j.cma.2005.08.026}}.
\newline\urlprefix\url{https://www.sciencedirect.com/science/article/pii/S0045782505005153}

\bibitem{RAPUN2010}
M.-L. Rapún, J.~M. Vega,
  \href{https://www.sciencedirect.com/science/article/pii/S0021999109007153}{Reduced
  order models based on local pod plus galerkin projection}, Journal of
  Computational Physics 229~(8) (2010) 3046--3063.
\newblock \href {https://doi.org/https://doi.org/10.1016/j.jcp.2009.12.029}
  {\path{doi:https://doi.org/10.1016/j.jcp.2009.12.029}}.
\newline\urlprefix\url{https://www.sciencedirect.com/science/article/pii/S0021999109007153}

\bibitem{2001219}
\href{https://www.sciencedirect.com/science/article/pii/S0950140101800458}{Chapter
  7 biot's theory for porous media}, in: J.~M. Carcione (Ed.), Wave Fields in
  Real Media: Wave Propagation in Anistropic, Anelastic and Porous Media,
  Vol.~31 of Handbook of Geophysical Exploration: Seismic Exploration,
  Pergamon, 2001, pp. 219--293.
\newblock \href {https://doi.org/https://doi.org/10.1016/S0950-1401(01)80045-8}
  {\path{doi:https://doi.org/10.1016/S0950-1401(01)80045-8}}.
\newline\urlprefix\url{https://www.sciencedirect.com/science/article/pii/S0950140101800458}

\end{thebibliography}

\end{document}